\documentclass[12pt]{article}
\usepackage{amsfonts}
\usepackage{amsmath}
\usepackage{amssymb}
\usepackage{latexsym}
\usepackage{amscd}
\usepackage{bbm}
\usepackage{url}

\title{On the interval of fluctuation of the singular values of random matrices
}
\author{
Olivier Gu\'edon${}^{1}$
\and
Alexander E. Litvak${}^{2}$
\and
Alain Pajor
       \and
Nicole Tomczak-Jaegermann${}^{3}$
}
\newcommand\address{\noindent\leavevmode

\medskip
\noindent
Olivier Gu\'{e}don, and Alain Pajor \\
Universit\'{e} Paris-Est \\
Laboratoire d'Analyse et Math\'{e}matiques Appliqu\'ees (UMR 8050) \\
UPEM, F-77454, Marne-la-Vall\'ee, France  \\
\texttt{\small olivier.guedon@u-pem.fr,  alain.pajor@u-pem.fr}

\medskip
\noindent
Alexander E. Litvak and Nicole Tomczak-Jaegermann, \\
Dept.~of Math.~and Stat.~Sciences,\\
University of Alberta, \\
Edmonton, Alberta, Canada, T6G 2G1.\\
\texttt{\small
e-mail:  aelitvak@gmail.com, nicole.tomczak@ualberta.ca}
}
\date{}

\newcounter{theorem}[section]
\newtheorem{theorem}{Theorem}
\newtheorem{lemma}[theorem]{Lemma}
\newtheorem{cor}[theorem]{Corollary}

\newtheorem{prop}[theorem]{Proposition}

\newcommand{\p}{\mathbb{P}}

\newcommand{\E}{\mathbb{E}}

\newcommand{\R}{\mathbb{R}}

\newcommand{\supp}{{\rm supp}}

\newcommand{\Rn}{\R^n}
\newcommand{\lc}{\lceil}
\newcommand{\rc}{\rceil}

\newcommand{\eps}{\varepsilon}

\newcommand{\lam}{\lambda}
\def\la{\left\langle}
\def\ra{\r\rangle}

\newcommand{\ta}{\tau}

\def\r{\right}

\newcommand{\qed}{\bigskip\hfill\(\Box\)}

\begin{document}
\maketitle

\footnotetext[1]{The research is partially
supported by the ANR GeMeCoD, ANR 2011 BS01 007 01}
\footnotetext[2]{Research partially supported by  the
E.W.R. Steacie Memorial Fellowship.}
\footnotetext[3]{This author holds the Canada Research Chair in
  Geometric Analysis.}

\begin{abstract}
  Let $A$ be a matrix whose columns $X_1,\dots, X_N$ are independent random
  vectors in $\R^n$. Assume that the tails of the 1-dimensional marginals
  decay as $\p(|\la X_i, a\ra|\geq t)\leq  t^{-p}$ uniformly in $a\in S^{n-1}$
  and $i\leq N$. Then for  $p>4$ we prove that with high probability $A/{\sqrt{n}}$
  has the Restricted Isometry Property (RIP) provided that Euclidean norms  $|X_i|$
  are concentrated around $\sqrt{n}$. We also show that the covariance matrix
  is well approximated by the empirical covariance matrix and 
  establish corresponding quantitative estimates on the 
 rate of convergence in terms of the ratio $n/N$. 
  Moreover, we obtain sharp bounds 
   for both problems when the decay is of the type $ \exp({-t^{\alpha}})$
  with $\alpha \in (0,2]$, extending the known case $\alpha\in[1, 2]$.
\end{abstract}



\section{Introduction and main results}
\label{intro}

Fix positive  integers $n, N$ and let  $A$ be an $n\times N$ random matrix
whose columns   $X_1,\dots, X_N$ are independent random  vectors in $\R^n$.
For a subset $I\subset\{1,\dots, N\}$ of cardinality $m$, denote by $A^I$
the $n\times m$ matrix whose columns  are $X_i,i\in I$. We are interested in
estimating the interval of fluctuation of the spectrum of some matrices related
to $A$ when the random vectors $X_i$, $i\leq  N$ have heavy tails; firstly,
uniform estimates of the spectrum of $(A^I)^\top A^I$ which is the set of
squares of the singular values of $A^I$, where $I$ runs over all subsets of
cardinality $m$  for some fixed parameter $m$ and secondly estimates for the
spectrum of $AA^\top$. The first problem is related to the notion of
Restricted Isometry Property (RIP) with $m$ a parameter of sparsity
whereas the second is about approximation of a covariance matrix by
empirical covariance matrices.

These questions have been substantially developed over recent years and
many papers devoted to these notions  were written. In this work, we say
that a random vector $X$ in $\R^n$ satisfies the tail behavior ${\bf H}(\phi)$ with parameter 
$\ta\geq 1$, if
 \begin{equation}
\label{condition-one}
  {\bf H}(\phi):\quad   \forall
  a \in S^{n-1}\ \forall t>0\quad \p\left( |\la X, a\ra |\geq t\r)
    \leq \ta/\phi(t)
\end{equation}
for a certain function $\phi$ and we assume that $X_i$ satisfies
$ {\bf H}(\phi)$ for all $i\leq N$.
We will focus on two choices of the  function $\phi$, namely
$\phi (t)=t^p$,  with  $p >4$, which means heavy tail behavior
for marginals, and $\phi (t)=(1/2) \exp(t^{\alpha})$,
with $\alpha \in (0,2]$, which corresponds to  an exponential
power type tail behavior and extends 
the known subexponential case ($\alpha=1$, see \cite{ALPT, ALPT1}).

The concept of the Restricted Isometry Property was introduced in
\cite{CT1} in order to study an exact reconstruction problem  by
$\ell_1$ minimization algorithm, classical in compressed sensing. Although it provided only a
sufficient condition for the 
reconstruction, it played a decisive role in the
development of the theory, and it is still an important property.
This is mostly due to the fact that a large number of important
classes of random matrices have RIP.
It is also noteworthy that the problem of reconstruction can be
reformulated in terms of convex geometry, namely in terms of
neighborliness of the symmetric convex hull of $X_1,\dots, X_N$, as
was shown in \cite{D3}.

Let us recall the intuition of RIP (for the definition see (\ref{rip-def})
below). For an $n\times N$ matrix $T$ and $1\le m \le N$, the {\it
isometry constant of order $m$} of $T$ is the parameter $0 < \delta_m(T)<1$
such that the square of Euclidean norms $|Tz| $ and $|z|$ are
approximately equal, up to a factor $1+\delta_m(T)$, for all $m$-sparse vectors
$z\in \R^N$ (that is, $|\supp (z)|\le m$). Equivalently, this means
that for every $I\subset \{1,\dots, N\}$ with $|I|\leq m$, the spectrum of
$(T^I)^\top T^I$ is contained in the interval
$[1-\delta_m(T),1+\delta_m(T)]$.  In particular when $\delta_m(T)<\theta$
for small $\theta$, then the squares of singular values of the matrices $T^I$
belong to $[1-\theta, 1+\theta]$. Note that in compressed sensing for the
reconstruction of vectors by $\ell_1$ minimization, one does not need RIP for
all $\theta>0$ (see \cite{D3} and \cite{lama}).
The RIP contains implicitly a normalization, in particular it implies that
the Euclidean norms of the columns  belong to an interval centered around one.

Let $A$ be an $n\times N$ random matrix  whose columns are $X_1, \ldots, X_N$.
In view of the example of  matrices with i.i.d. entries, centered and with variance
one, for which $\E|X_i|^2=n$, we  normalized the matrix  by considering $A/\sqrt n$
and we  introduce the concentration function
\begin{equation} \label{conc}
      P(\theta):= \p \left( \max_{i\le
      N}\left|{|X_i|^2\over n}-1\right|\geq \theta \right) .
\end{equation}

Until now the only known cases of random matrices satisfying a RIP were the cases of
subgaussian \cite{CRT, CT1, D3, MPT} and  subexponential \cite{ALPT2} matrices.
Our first main theorem says that matrices we consider
have the RIP of order $m$, with ``large" $m$ of the form $m= n\psi(n/N)$
with $\psi$ depending on $\phi$ and possibly on other parameters. In particular,
when $N$ is proportional to $n$, then $m$ is proportional to $n$.
We present a simplified version of our result, for the
detailed version see Theorem~\ref{thm:RIP-M} below.

\begin{theorem}
\label{RIP-intro}
Let $0 < \theta < 1$. Let $A$ be  a random  $n\times N$ matrix  whose columns
$X_1, \ldots, X_N$ are independent random vectors satisfying hypothesis $\mbox{\bf H}(\phi)$
for some $\phi$.  Assume that $n, N$ are large enough.
Then there exists a function $\psi$ depending on $\phi$ and $\theta$ such that
with high probability (depending on the concentration function $P(\theta)$)
the matrix $A/\sqrt{n}$ has RIP of order $m =[ n\psi(n/N)]$ with a
parameter $\theta$ (that is,
$\delta _m(A/\sqrt{n})\leq \theta$).
\end{theorem}

The second problem we investigate goes back to a question of
Kannan, Lov{\'a}sz and Simonovits (KLS).
As before assume that $A$ is a random  $n\times N$ matrix  whose columns 
$X_1, \ldots, X_N$ are independent 
random vectors satisfying hypothesis $\mbox{\bf H}(\phi)$ for some $\phi$. 
Additionally assume that $X_i$'s are identically distributed as 
a centered random vector $X$.  KLS
question asks how fast the empirical covariance matrix $U:=(1/N)
AA^\top$ converges to the covariance matrix $\Sigma := (1/N) \E
AA^\top= \E U$. Of course this depends on assumptions on $X$.
In particular, is it true that with high probability the
operator norm $\|U - \Sigma\| \leq \eps \|\Sigma\|$ for $N$ being
proportional to $n$?  Originally this was asked  for  log-concave
random vectors but the general question of approximating the
covariance matrix by sample covariance matrices is an important
subject in Statistics as well as on its own right. The corresponding  question
in Random Matrix Theory is about the limit behavior of smallest and
largest singular values. In the case of Wishart matrices, that is
when the coordinates of $X$ are i.i.d. centered random variables of
variance one, the Bai-Yin theorem \cite{BY} states that under assumption
of boundedness of fourth moments the limits
of minimal and maximal singular numbers of $U$ are $(1\pm\sqrt{\beta})^2$ as
$n,N\to\infty$ and $n/N\to \beta\in(0,1)$. Moreover, it is known
\cite{BSY, Silv} that boundedness of fourth moment is necessary in order to
have the convergence of the largest singular value. The asymptotic
non-limit behavior (also called ``non-asymptotic'' in Statistics),
i.e., sharp upper and lower bounds for singular values in terms of
$n$ and $N$, when $n$ and $N$ are sufficiently large, was studied in
several works.  To keep the notation more compact and clear we put
\begin{equation}
  \label{eq:ms_def}
  M:= \max _{i\leq N} |X_i|,
\qquad S: =  \sup_{a\in S^{n-1}} \left| \frac{1}{N}
  \sum_{i=1} ^N\left(\langle X_i,  a\rangle^2 -  \E\langle X_i,
    a\rangle^2\right)\right|.
\end{equation}
Note that if $\E \langle X , a\rangle^2 =1$ for every $a\in S^{n-1}$
(that is, $X$ is isotropic), then the bound $S\leq \eps$ is equivalent
to the fact that the singular values of $U$
belong to the interval $[1-\varepsilon, 1+\varepsilon]$.
For Gaussian matrices it is known (\cite{edel,szarek})
that with probability close to one
\begin{equation}
\label{Bai-Yin quantitative}
 S \leq C \sqrt{n/N},
\end{equation}
where $C$ is a positive absolute constant.
In \cite{ALPT, ALPT1}  the same estimate was obtained for a
large class of random matrices,
which in particular did not require that entries of the columns are
independent, or that
$X_i$'s are identically distributed. In particular this  solved the
original KLS problem. More precisely,
(\ref{Bai-Yin quantitative}) holds with high probability
under the assumptions that the  $X_i$'s satisfy
hypothesis $\mbox{\bf H}(\phi)$
with $\phi (t)=e^t/2$ and that $M \leq C (Nn)^{1/4}$ with
high probability. Both conditions hold for log-concave random vectors.

Until recent time,
quite strong conditions on the tail behavior of the one dimensional
marginals of the $X_i$ were imposed, typically of subexponential type.
Of course, in view of Bai-Yin theorem, it is a natural question
whether one can replace the function $\phi (t)=e^t/2$ by the function
$\phi (t)=e^{t^{\alpha}}/2$ with $\alpha\in (0,1)$ or $\phi(t) = t^p$,
for $p\geq 4$.
The first attempt in this direction was done in \cite{V}, where the
bound $S\leq C(p, K) (n/N)^{1/2-2/p} (\ln \ln n)^2$ was obtained for
every $p>4$ provided that $M\leq K \sqrt{n}$. Clearly, $\ln \ln n$ is
a ``parasitic" term, which, in particular, does not allow to solve the
KLS problem with $N$ proportional to $n$.  This problem was solved in
\cite{MP, SV}  under strong assumptions  and in particular when
$M\leq K \sqrt{n}$ and $X$ has i.i.d. coordinates with bounded $p$-th moment
with $p>4$.
Very recently, in \cite{MPaouris2012}, the ``right"
upper bound $S\leq C (n/N)^{1/2}$ was proved for $p>8$ provided that
$M \leq C (Nn)^{1/4}$. The methods used in
\cite{MPaouris2012} play an influential role in the present paper.

The problems of estimating the smallest and the largest singular values
are quite different. One expects weaker assumption for estimating the smallest
singular value. This already appeared in the work \cite{SV} and was pushed 
further in recent works \cite{KM, Ti, Y} and in 
\cite{GeMu, LM, Ol} which led to new bounds 
on the performance of $\ell_1$-minimization methods.

In this paper we solve the KLS problem for $4<p\leq 8$, in Theorem~\ref{thm:KLS-M}.
Our argument works also in other cases  and makes the bridge between the known cases
$p>8$ and the  exponential case.

\begin{theorem}
  \label{thm:KLS-M}
Let $X_1,\dots, X_N$ be independent random  vectors in $\R^n$ satisfying
hypothesis $\mbox{\bf H}(\phi)$ with $\phi (t)=t^p$ for some $p\in (4, 8]$.
  Let $\eps \in (0, 1)$ and $\gamma = p-4-2\eps >0$.  Then
\begin{equation}\label{eqshs}
S \leq    C \left(\left(\frac{M^2}{n}\r) \left( \frac{n}{N}\r)
   + C(p, \eps)\,
    \left(\frac{n}{N}\r)^{\gamma/p}  \r),
\end{equation}
  with  probability larger than
$ 1- 8 e^{-n} -  2\eps^{-p/2} \max\{N^{-3/2}, n^{-(p/4-1)}\}$.
\end{theorem}

In particular, if $N$ is proportional to $n$ and $M^2/n$ is bounded by a constant with high probability, which is the case for  large classes of random vectors, then with high probability
$$
   S\leq C \left(n/N\right)^{\gamma/p}.
$$

Let $X$ have i.i.d. coordinates distributed as a centered random variable
with finite $p$-th moment, $p>2$. Then by Rosenthal's inequality (\cite{Ros}, see
also \cite{FHJSZ} and Lemma~\ref{rosin-l} below), 
$X$ satisfies hypothesis $\mbox{\bf H}(\phi)$ with
$\phi (t)=t^p$. Let $X_1$, ..., $X_N$ be independent random vectors distributed
as $X$. It is known (\cite{BSY}, \cite{Silv}, see also \cite{LS} for a quantitative
version) that when $N$ is proportional to $n$ and in the absence of fourth moment,
$M^2/n\to\infty$ as $n\to \infty$. Hence, bounds for $S$ involving the term $M^2/n$  
like the bound (\ref{eqshs}) are of interest 
only for $p\geq 4$. We don't know if it holds in the case $p=4$.

The main novelty of our proof is a delicate analysis of the behavior of norms of
submatrices, namely quantities $A_k$ and $B_k$, $k\leq N$, defined in (\ref{eq:AkBk_intro})
below. This analysis is done in Theorem~\ref{A_k}, which is  in the heart of the
technical  part of the paper and it will be presented in the next section.
The estimates for $B_k$ are responsible for RIP, Theorem~\ref{RIP-intro}, while
the estimates for $A_k$ are responsible for KLS problem,  Theorem~\ref{thm:KLS-M}.

\medskip

As usual in this paper $C$, $C_0$, $C_1$, ..., $c$, $c_0$, $c_1$, ...  always
denote absolute positive constants whose values may vary from line to line.

\vspace{2ex}

The paper is organized as follows. In Section~\ref{main-tech}, we formulate
the main technical result. For the reader convenience, we postpone its proof
till Section~\ref{sec:proof_C}. In Section~\ref{rip-M}, we discuss the results
on RIP. The fully detailed  formulation of the main result in this direction
is Theorem~\ref{thm:RIP-M}, while Theorem~\ref{RIP-intro} is its very simplified
corollary. In Section~\ref{klm-M}, we prove Theorem~\ref{thm:KLS-M} as a consequence
of Theorem \ref{thm:KLS}.  The case $p>8$ and the exponential cases are proved in
Theorem \ref{thm:KLS-new} using the same argument. Symmetrization and formulas for
sums of the $k$ smallest order statistics of independent non-negative random variables
with heavy tails allow to reduce the problem on hand to estimates for $A_k$.
In the last Section~\ref{sec:optimal}, we discuss optimality of the results.

\smallskip

An earlier version of the main results of this paper was announced in \cite{GLPT}.

\bigskip

\noindent {\bf Acknowledgment.}
A part of this research was performed while the authors were visiting at several
universities.  Namely, the first named author visited University of Alberta at
Edmonton in April 2013 and the second and the fourth named author visited University
Paris-Est  in June 2013 and in June 2014. The authors would like to thank
these universities for their support and hospitality.

\section{Norms of submatrices}
\label{main-tech}

We start with a few general preliminaries and notations.
We denote by $B_2^n$ and $S^{n-1}$  the standard unit Euclidean ball and the
unit sphere in $\R^n$
and by $|\cdot|$  and $\la\cdot, \cdot \ra$  the corresponding
Euclidean norm and inner product. Given a set $E\subset \{1, ..., N\}$,
$|E|$ denotes its cardinality and $B_2^E$ denotes the unit Euclidean ball
in $\R^E$, with the convention $B_2^{\emptyset}=\{0\}$.

A standard volume argument implies that
for every integer $n$ and for every $\eps \in(0,1)$ there exists an $\eps$-net
$\Lambda\subset B_2^n$ of $ B_2^n$
of cardinality not exceeding $(1 + 2/\eps)^n$; that is, for every $x\in B_2^n$,
$\min_{y\in\Lambda}|x-y|<\varepsilon$. In particular, if $\eps \leq 1/2$
then the cardinality of $\Lambda$ is not larger than $(2.5/\eps)^n$.

\vspace{2ex}

By $\cal{M}$ we denote the class of increasing functions
$\phi : [0, \infty) \to [0, \infty)$ such that the function
$\ln \phi \left(1/\sqrt{x}\r)$  is convex on $(0,\infty)$.
The  examples of such functions  considered in this paper  are
$\phi(x)=x^p$ for some $p>0$ and $\phi(x)=(1/2)\exp(x^\alpha)$ for
some $\alpha>0$.

Recall that the  hypothesis ${\bf H}(\phi)$ has been  defined in
the introduction by (\ref{condition-one}). Note that this
hypothesis is satisfied if
$$
    \sup _{a\in S^{n-1}}\E\, \phi(|\la X, a\ra |) \leq \ta.
$$

For $k\leq N$ and random vectors $X_1$, ..., $X_N$ in $\R^n$ we
define $A_k$ and $B_k$ by
\begin{equation}
  \label{eq:AkBk_intro}
     A_k:= \sup_{a\in S^{N-1}\atop {|\supp (a)|\le k}} \left|
     \sum _{i=1}^{N} a_i X_i \r|, \quad B_k^2 := \sup_{a\in
     S^{N-1}\atop {|\supp (a)|\le k}} \left|\left| \sum _{i=1}^{N}
     a_i X_i \r|^2 - \sum _{i=1}^{N} a_i^2 |X_i|^2\r|.
\end{equation} 
We would like to note that $A_k$ is the supremum of norms of submatrices
consisting of $k$ columns of $A$, while $B_k$ plays a crucial role
for RIP estimates. We provide more details on the role of $A_k$ and 
$B_k$ in the next section.

Recall also a notation from the introduction
$$
 M=\max_{i\leq N} |X_i|.
 $$

\vspace{1ex}

We formulate now the main technical result, Theorem~\ref{A_k},
which is the key result for both bounds for $A_k$ and for $B_k$.
The role of $A_k$ and $B_k$ in RIP estimates will be explained
in the next section. We postpone the proof to Section~\ref{sec:proof_C}.

\begin{theorem}
       \label{thm:A_k-M} \label{A_k}
Let  $p>4$, $\sigma \in (2, p/2)$, $\alpha \in (0,2]$, $t>0$, and
$\ta, \lambda \geq 1$. Let $X_1,\dots, X_N$ be independent random  
vectors in $\R^n$ satisfying hypothesis $\mbox{\bf H}(\phi)$ with parameter $\ta$
 either for $\phi (x) =x^p$ or for $ \phi (x) = (1/2) \exp (x^{\alpha})$. For $k\leq N$
define $M_1$, $\beta$ and $C_{\phi}$ in two following cases.

\noindent
{\bf Case 1. } $\phi (x) = x^p$.
We   assume that $\lambda \leq p$ and we let $C_{\phi} = e^4$,
$$
   M_1 :=
   C_1(\sigma,\lambda, p)
  \sqrt{k} \, \left(\frac{N\ta }{k}\r)^{\sigma/p}
\quad \mbox{ and } \quad
  \beta :=   C_2(\sigma, \lambda) (\ta N)^{-\lambda}
     +  C_3(\sigma, \lambda, p)
     \frac{N^2 \ta}{ t ^p},
$$
where
$$
   C_1(\sigma, \lambda, p) =
   32\, e^4\, \sqrt{\frac{\sigma+\lambda}{1+\lambda/2}}\,
  \left(\frac{2p}{p-2\sigma}\right)^{1+2\sigma /p}
  \left(\frac{\sigma+\lambda}{\sigma-2}\right)^{2\sigma/p} \,
   (20e)^{\sigma/p},
$$
$$
  C_2(\sigma, \lambda) :=   \left(\frac{2(\sigma+\lambda)}{5e (\sigma-2)}\right)^{\lambda}
  \frac{1}{2\lambda -1}  \quad \mbox{and} \quad
  C_3(\sigma, \lambda, p) :=  \frac{(\sigma +\lambda)^p}{4  (2 (\sigma -2))^p} .
$$

\noindent
{\bf Case 2. } $ \phi (x) = (1/2) \exp (x^{\alpha})$.
We assume that $\lambda \geq 2$ and  we let $C_{\phi} = C^{1/\alpha}$,
where $C$ is an absolute positive constant,
$$
   M_1 :=  (C\lambda)^{1/\alpha}\,
 \sqrt{ k}\,
    \left( \ln \frac{2 N\ta}{ k} + \frac{1}{\alpha}\right)^{1/ \alpha}
$$
and
$$
 \beta:= \frac{1}{(10 N\ta )^{\lambda}}  \exp\left(- \frac{\lambda k^{\alpha/2}}{(3.5\,
  \ln(2 k))^{2\alpha}}\r)  + \frac{N^2\ta }{2\exp((2t)^{\alpha})} .
$$
 In both cases we also assume that $\beta <1/32$.
Then with probability at least $1- \sqrt\beta$ one has
$$
  A_k \leq (1-4\sqrt{\beta})^{-1}\left(M + 2  \sqrt{C_{\phi}\, t\, M}
            + M_1 \r)
$$
and
\begin{align*}
     B_k^2 \leq  (1-4\sqrt{\beta})^{-2}\left(4\sqrt{\beta}
     M^2 +  \left(8 C_{\phi}\, t + M_1\right) M + 2 M_1^2 \r).
\end{align*}
\end{theorem}

We would like to emphasize that  $A_k$ and $B_k$ are of
different nature. In particular, Theorem~\ref{A_k} in the case $\phi (x)=x^p$
has to be applied with different choices of the parameter $\sigma$. We summarize those
choices in the following remark.

\medskip

\noindent
{\bf Remark. } In the case $\phi(x)=x^p$ we will use the following two
choices for $\sigma$: \\
 {\bf 1.}  Choosing $\sigma =p/4$ and
assuming $p>8$ we get
$$
  M_1 \leq   C \sqrt{\frac{p}{\lambda}}\,
  \sqrt{\frac{p}{p-8}} \, \sqrt{k} \, \left(\frac{N\ta }{k}\r)^{1/4}
$$
and
$$
  \beta \leq  \left(\frac{2 (p+4\lambda)}{5 e N \ta (p-8)}\right)^{\lambda}
  \frac{1}{2\lambda -1}  + \frac{N^2 \ta (p+4\lam)^p}{4  (2t (p -8))^p}.
$$
{\bf 2.}
Choosing $\sigma = 2+ \eps$ with $\eps\leq \min\{1, (p-4)/4\}$,  we get
$$
  M_1 \leq C \,
  \left(\frac{p}{p-4}\right)^{1+(4+2\eps) /p}
  \left(\frac{\lambda}{\eps}\right)^{2(2+\eps)/p} \, \sqrt{k} \, 
  \left(\frac{N\ta }{k}\r)^{(2+\eps)/p}
$$
and
$$
  \beta =   \left(\frac{2(3+\lambda)}{5e \eps N\ta }\right)^{\lambda}
  \frac{1}{2\lambda -1}  + \frac{N^2 \ta (3 +\lambda)^p}{4  (2\eps t)^p}.
$$

\medskip

\noindent
{\bf Remarks on optimality.} \\
{\bf 1.  The case} $\phi(x)=x^p$, $p>4$. Let $\ta\geq 1$, 
$N\geq (64 C_2(\sigma, \lambda))^{1/\lambda}$ and 
$t= (64 N^2 C_3(\sigma, \lambda, p))^{1/p}$.   
Then $\beta \leq 1/32$ and $\sqrt{tM}\leq C_4(\sigma, \lambda, p) (M+M_1)$. 
Hence with probability larger than 3/4 we have
$$
  A_k \leq  C(\sigma, \lambda, p)
  \left(M + \sqrt{k} (N\ta /k)^{\sigma/p} \r).
$$
In Proposition~\ref{6.5} below we show that there exist independent random vectors $X_i$'s 
satisfying the conditions of Theorem~\ref{A_k} with $\ta=1$ and such that
$$
  A_k \geq C(p) \sqrt{k} (N/k)^{1/p} \left(\ln (2N/k) \r) ^{-1/p}
$$
with probability at least 1/2. 
Note that $M=A_1\leq A_k$.  Therefore 
\begin{equation}\label{sharp1}
\max\{ M, C(p) \sqrt{k} (N/k)^{1/p} \left(\ln (2N/k) \r) ^{-1/p}\} \leq
A_k \leq  C(\sigma, \lambda, p) \left(M +\sqrt{k} (N/k)^{\sigma/p}\r)
\end{equation}
with probability at least 1/4.

\medskip

\noindent
{\bf 2. The case} $ \phi (x) = (1/2) \exp (x^{\alpha})$, $\alpha\in[1,2]$. Let
$\lambda =2$ and $t= (\ln N)^{1/\alpha}$. Then $\beta\leq 1/32$. 
Hence with probability larger than 3/4 we have
$$
  A_k \leq C\left(M + C^{1/\alpha} \sqrt{k} \left(\ln (6^{1/\alpha}\ta N/k)\r)^{1/\alpha} \r).
$$
In Proposition~\ref{6.7} below we show that there exist independent random vectors $X_i$'s 
satisfying the conditions of Theorem~\ref{A_k} with $\ta$ bounded by an absolute constant 
and such that
$$
  A_k \geq  \sqrt{k/2}  \left(\ln (N/(k+1)) \r) ^{1/\alpha}
$$
with probability at least 1/2. Using again that  $M=A_1\leq A_k$ we observe 
\begin{equation}\label{sharp2}
\max\{M,  \sqrt{k/2}  \left(\ln (N/(k+1)) \r) ^{1/\alpha}\} \leq
   A_k \leq C\left(M + C^{1/\alpha} \sqrt{k} \left(\ln (6^{1/\alpha} N/k)\r)^{1/\alpha}\r)
\end{equation}
with probability at least 1/4.

\section{Restricted Isometry Property}
\label{rip-M}

We need more definitions and notations.

Let $T$ be an $n\times N$ matrix and let
$1\le m \le N$.  The $m$-th {isometry constant} of $T$ is defined as the smallest number $\delta_m=\delta_m(T)$ so that
\begin{equation}
  \label{rip-def}
  (1-\delta_m)|z|^2\le {|T z|^2}    \le( 1+\delta_m) |z|^2
\end{equation}
holds for all vectors $z\in\R^N$ with $|\supp (z)|\le m$.
For $m=0$, we put $\delta_0(T)=0$. Let $\delta\in (0,1)$.
The matrix $T$ is said to satisfy the Restricted Isometry
Property of order $m$ with parameter $\delta$, in short
$\ {\rm RIP}_m (\delta)$, if $0\leq\delta_m(T)\leq \delta$.

Recall that a vector $z \in \R^N$ is called $m$-sparse if $|\supp (z)|\le m$.
The subset of $m$-sparse unit vectors in $\R^N$ is denoted by
$$
   U_m=U_{m}(\R^N): =\{z\in\R^N\,:\, |z|=1, | \supp (z)|\le m\}.
$$

Let $X_1$, ..., $X_N$ be random vectors in $\R^n$ and let $A$ be the
$n\times N$ matrix whose columns are the $X_i$'s.
By the definition of $B_m$ (see (\ref{eq:AkBk_intro}))
we clearly have
$$
\max_{i \le N} \left|{|X_i|^2\over n} -1\right| =
\delta_1\left({A\over\sqrt n}\right) \leq \delta_m\left({A\over\sqrt
    n}\right) = \sup_{z\in U_m}\left| {|Az|^2\over n}-1 \right|
$$
\begin{equation}\label{RIPfromBm}
 \le {B_m^2\over n}+ \max_{i \le  N} \left|{|X_i|^2\over n} -1\right|.
\end{equation}
Thus, in order to have a good bound on $\delta_m\left({A/\sqrt
    n}\right)$ we require a strong concentration of each $|X_i|$ around
$\sqrt n$ and we need to estimate $B_m$.

To control the concentration  of $|X_i|$ we consider the function $P(\theta)$,
defined in the introduction by (\ref{conc}). 
Note that this function estimates the concentration of the maximum. Therefore,
when it is small, we have much better concentration of each $|X_i|$ around
$\sqrt{n}$.

We are now ready to state the main result about RIP.
Theorem~\ref{RIP-intro}, announced in the introduction, is
a very simplified form of it.

\begin{theorem}
\label{thm:RIP-M}
Let  $p>4$,  $\alpha \in (0,2]$, $\ta\geq 1$ and $1\leq n\leq N$.
Let $X_1,\dots, X_N$ be independent random  vectors in $\R^n$ satisfying 
hypothesis $\mbox{\bf H}(\phi)$ with the parameter $\ta$ either for 
$\phi (x) =x^p$ or for $ \phi (x) = (1/2) \exp (x^{\alpha})$. 
Let  $P(\cdot)$ be as in  (\ref{conc}) and  $\theta \in (0, 1)$.

\noindent
{\bf {Case 1.}}
  $\phi(x)=x^p$. Let $\eps \leq \min\{1,
 (p-4)/4\}$. Assume that
$$
  \frac{2^8}{\eps\, \theta\, \ta} \leq N
 \leq c\, \theta  \left(c\, \eps\, \theta \r)^{p/2}\, n^{p/4} \sqrt{\ta}
$$
and set
$$
     m = \left[ C(\theta, \eps, p) \, n \left(\frac{N\ta }{n}\right)
       ^{-2(2+\eps)/(p-4-2\eps)} \r]   \, \, \mbox{ and } \, \,
    \beta = \frac{4}{3 e^2 \, \eps^2\, N^2 \ta^2} + \frac{5^p \, N^2 \ta}{4
    (2 c \, \eps\, \theta )^p\, n^{p/2}}   ,
$$
where
\begin{equation}
\label{eq:def_C}
     C(\theta, \eps, p) =  c \left(\frac{p-4}{p}
       \r)^{2(p+4+2\eps)/(p-4-2\eps)} \,
     \eps ^{4(2+\eps)/(p-4-2\eps)}  \, \theta ^{2p/(p-4-2\eps)} ,
\end{equation}
$c$ and $C$ are  absolute positive constants.

\smallskip

\noindent
{\bf {Case 2.}} $ \phi (x) = (1/2) \exp (x^{\alpha})$.
Assume that
$$
 \frac{1}{\ta}\, \max\left\{2^{1/\alpha}, 4/\theta \r\} \leq N
 \leq  c\, \theta\, \sqrt{\ta}\,   \exp\left((1/2)\left( c\,
    \theta     \sqrt{n}\r)^{\alpha}\r)
$$
and set
$$
    m=\left[ C^{-2/\alpha}\, \theta ^2 \, n \left(\ln (C^{2/\alpha}\,
        N\ta / (\theta ^2 \, n))\r)^{-2/\alpha} \r]
$$
and 
$$
   \beta = \frac{1}{(10 N \ta)^2 }\exp{\left(\frac{-2 m^{\alpha/2}}{(3.5 \ln(2m))^{2\alpha}}\r)}
   + \frac{N^2\ta }{2} \exp\left(- c \left(\theta \sqrt{n}\r)^{\alpha}\r)   ,
$$
where $c$ and $C$ are  absolute positive constants.

\noindent
Then  in both cases we have 
$$
   \p\left( \delta _m(A/\sqrt{n})\leq \theta\r) \geq 1 - \sqrt \beta
        - P(\theta/2) .
$$
\end{theorem}

\noindent {\bf Remarks. 1.}
Note that  for instance in case 1, the constraint  $N\leq c(\theta, \eps, \ta, p) n^{p/4}$
is not important because for $N\gg n^{p/4}$ one has
$$
  m = \left[ C(\theta, \eps, p) \, n \left(\frac{N\ta }{n}\right)
       ^{-2(2+\eps)/(p-4-2\eps)} \r] = 0.
$$
 A similar remark is valid in the second case.
\newline
{\bf 2. }
In most applications $P(\theta)\to 0$ very fast as $n, N\to \infty$.
For example, for so-called isotropic log-concave random vectors it
follows from results of Paouris (\cite{Pao, Paor2}, see also \cite{Klartag,
GM} or Lemma~3.3 of \cite{ALPT2}).
As another example consider the model when $X_i$'s are i.i.d. and moreover
the coordinates of $X_1$ are i.i.d. random variables distributed as a random
variable $\xi$. In the case when $\xi$ is of variance one and has finite
$p$-th moment, $p>4$, then by Rosenthal's inequality $P(\theta)$ is well bounded
(for a precise bound see Corollary~\ref{ros-conc} below,
see also Proposition~1.3 of \cite{SV}). Another case is when $\xi$ is the Weibull
random variable of variance one, that is consider $\xi _0$ such that
$\p\left(|\xi_0 | > t\r)= \exp{(- t^{\alpha})}$ for $\alpha\in (0, 2]$ and
let $\xi =\xi_0/ \sqrt{\E\xi_0^2}$. By Lemma~3.4 from \cite{ALPT2}
(see also Theorem~1.2.8 in \cite{lama}),  $P(\theta)$ satisfies
(\ref{condtwoconc}) below.
\newline
{\bf 3. Optimality.}  Taking $\eps$ in the Case 1 of order $(p-4)^2/\ln (2N\ta/n)$ 
and assuming that it satisfies the condition of the theorem, we observe that in Case~1
$$
   m = \left[ C(\theta, p) \, n \left(\frac{N\ta}{n}\right)
       ^{-4/(p-4)} \, \left(\ln \frac{2N\ta }{n}\r)^{-8/(p-4)} \r] .
$$
Moreover, Proposition~\ref{6.6} below shows that for $q>p>4$ there are independent 
random vectors $X_i$'s satisfying hypothesis $\mbox{\bf H}(\phi)$ with parameter 
$\ta=\ta(p, q)$ and such that for $\theta= 1/2$, 
$N \leq C(p, q) n^{p/4}\left(\ln (2N/n) \r)^{-p/2}$ one can't get better estimate than
$$
 m\leq 8 \left(N/n\r)^{-2/(q-2)} n .
$$
{\bf 4. Optimality.} In Case 2 with $\ta$ bounded by an absolute constant and $\alpha\in[1,2]$, 
let $c_0 n \leq N \leq c_1 \exp(c_2 n^{\alpha/2})$, $\theta = 0.4$ and assume that $P(\theta/2)$ 
is small enough. Then $\p (\delta _m \leq 1/2)\geq 1/2$ provided that 
$m=n\left(C\ln(C^{2/\alpha} N/n)\r)^{2/\alpha}$. Proposition~\ref{6.7} below shows that 
the estimate for $m$ is sharp, that is in general $m$ can't be larger than 
$m=n \left(C\ln(2 N/n)\r)^{2/\alpha}$.

\medskip

\noindent
{\bf Proof. }
We first pass to the subset $\Omega _0$ of our initial probability
space where
$$
   \max_{i\le N} \left| {|X_i|^2\over n}-1 \right|\leq \theta/2.
$$
Note that by (\ref{conc}) the probability of this event is at least
$1-P(\theta/2)$ and if this event occurs then we also have
$$
  \max_{i\le N}{|X_i|} \leq  3\sqrt{n}/2 .
$$

We will apply Theorem~\ref{A_k} with $k=m$, $t=  \theta \sqrt{n}/(100 C_{\phi})$,
where $C_{\phi}$ is the constant from Theorem~\ref{A_k}. Additionally we assume that
$\beta \leq 2^{-9}\theta ^2$ and $M_1\leq t$.  Then with
probability at least $1-\sqrt \beta -P(\theta /2)$ we have
$$
   B_m^2 \leq (16 \sqrt{\beta}\,  + \theta /4 ) n \leq \theta \, n / 2 .
$$
Together with (\ref{RIPfromBm}) this proves $\delta _m(A/\sqrt{n})\leq \theta$.
Thus we only need to check when the estimates for $\beta$ and
$M_1$ are satisfied.

\smallskip

\noindent
{\bf Case 1. $\phi(x)=x^p$. } We start by proving the estimate for $M_1$. 
We let
$\sigma = 2+ \eps$,  $\eps \leq \min\{1, (p-4)/4\}$ and  $\lambda
=2$. Then by  Theorem~\ref{A_k}
(see also the Remark following it),
for some absolute constant $C$ we have
$$
     M_1\leq C \,
  \left(\frac{p}{p-4}\right)^{1+(4+2\eps) /p}
  \left(\frac{1}{\eps}\right)^{2(2+\eps)/p} \, \sqrt{m} \, \left(\frac{N\ta}{m}\r)^{(2+\eps)/p} .
$$
Therefore the estimate   $M_1\leq c \, \theta \sqrt{n}$ with $c= 1/(100 e^4)$ is satisfied
provided that
$$
   m = \left[ C(\theta, \eps, p) \, n \left(\frac{N\ta}{n}\right)
     ^{-2(2+\eps)/(p-4-2\eps)} \r] ,
$$
with $C(\theta, \eps, p)$ defined in (\ref{eq:def_C}) and
the absolute constants  properly adjusted.

Now we estimate the  probability. From Theorem~\ref{A_k}
(and the Remark following it),
with our choice of $t$ and $\lambda$ we have
$$
  \beta \leq \frac{4}{3 e^2 \, \eps^2\, N^2\, \ta^2 } + \frac{5^p \, N^2\ta  }{4
    (2 c \, \eps\, \theta )^p\, n^{p/2}}
  \leq 2^{-9} \theta^2
$$
provided that $2^8/(\eps\, \theta\ta) \leq N \leq 2^{-4}\theta\sqrt{\ta} \left(0.4
  c \, \eps\, \theta \r)^{p/2}\, n^{p/4}$.  This completes the proof
  of the first case.

\medskip

\noindent
{\bf Case 2. $\phi(x)=(1/2)\exp(x^\alpha)$.\, }
As in the first case we start with the condition
$M_1\leq t$.
We choose $\lambda =4$. Note that  $N\ta /m\geq 2^{1/\alpha}$ as
$N\ta \geq 2^{1/\alpha} n$. Therefore for some absolute constant $C$,
$$
     M_1\leq \sqrt{m} \left( C \ln (2N\ta/m) \r)^{1/\alpha}.
$$
Therefore the condition  $M_1\leq t$ is satisfied
provided that
$$
     m\leq C_1^{-2/\alpha}\, \theta ^2 \, n \left(\ln
       (C_1^{2/\alpha}\, N\ta /(\theta ^2 \, n))\r)^{-2/\alpha}
$$
for an absolute positive constant $C_1$. This justifies the choice of
$m$.

Now we estimate  the  probability.  From Theorem~\ref{A_k}
with our choice of $t$ and $\lambda$ we have
$$
   \beta \leq \frac{1}{(10 N\ta)^2 }
   \exp{\left(\frac{-2 m^{\alpha/2}}{(3.5 \ln(2m))^{2\alpha}}\r)}
   + \frac{N^2\ta }{2} \exp\left(- c \left(\theta \sqrt{n}\r)^{\alpha}\r)
  \leq 2^{-9} \theta^2 ,
$$
provided that $4/(\theta \ta) \leq N \leq 2^{-5}\theta \sqrt{\ta}
\exp\left( c \left(\theta \sqrt{n}\r)^{\alpha}\r) $.
This completes the proof.
\qed

\section{Approximating  the covariance matrix}
\label{klm-M}

We start with the following $\eps$-net argument for bilinear forms, 
which will be used below.

\begin{lemma} \label{matrix}
Let $m\ge 1$ be an integer and $T$ be an  $m\times m$ matrix.  Let
$\eps \in (0,1/2)$ and
${\cal{N}}\subset B_2^m$ be an $\eps$-net  of $B_2^m$ (in the Euclidean
metric).
Then 
$$
 \sup_{x\in B_2^m}|\langle T x,x\rangle|\leq (1-2\eps)^{-1}\,
 \sup_{y\in {\cal{N}}}|\langle T y,y\rangle|.
$$
\end{lemma}

\noindent
{\bf Proof. }
Let $S=T+T^*$. For any $x,y\in\R^m$,
$$
 \langle Sx,x\rangle=\langle Sy,y\rangle+\langle Sx,x-y\rangle
+\langle S(x-y),y\rangle .
$$
Therefore
$|\langle Sx,x\rangle|\leq |\langle Sy,y\rangle|+2|x-y| \|S\|$.
Since $S$ is symmetric, we have
$$
  \|S\|=\sup_{x\in B_2^m}|\langle Sx,x\rangle|.
$$
Thus, if $|x-y| \leq \eps$,
then
$$
  \|S\|\leq \sup_{y\in {\cal{N}}}|\langle Sy,y\rangle| +2\eps \|S\|
$$
and
$$
  \sup_{x\in B_2^m}|\langle Sx,x\rangle|\leq (1-2\eps)^{-1}\,
  \sup_{y\in {\cal{N}}}|\langle Sy,y\rangle|.
$$
Since $T$ is a real matrix, then for every $x\in\R^m$,
$\langle Sx,x\rangle=2\langle T x,x\rangle$. This concludes the proof.
\qed

Now we can prove the following technical lemma, which 
emphasizes the role of  the parameter $A_k$ in estimates of the
distance between the covariance matrix and the empirical one.  This
role was first recognized in \cite{Bourgain} and \cite{ALPT}.
Other  versions of the lemma appeared in \cite{ALPT2, ALLPT}.
Its proof uses the symmetrization method as in \cite{MPaouris2012}.

\begin{lemma}
\label{KLS technic}
Let $\ta\geq 1$, $1\leq k< N$ and $X_1,\dots, X_N$ be  independent random vectors in $\R^n$.
Let  $p\geq 2$, $\alpha \in (0,2]$. Let $\phi$ be either $\phi (t)=t^p$
in which case we set $C_{\phi} = 8 \ta^{2/p} N^{2/\min(p,4)}$ and assume
$$
  \forall 1\leq i\leq N\quad\forall a\in S^{n-1}\quad
   \E \,  |\langle X_i, a\rangle|^p\leq \ta,
$$
or $\phi (t)=(1/2) \exp(t^{\alpha})$ in which case we assume that
$X_i$'s satisfy hypothesis $\mbox{\bf H}(\phi)$ with parameter $\ta$ and set
$C_{\phi} = 8 \sqrt{C_{\alpha} N\ta}$, where
$C_{\alpha}=(8/\alpha) \, \Gamma (4/\alpha)$, $\Gamma (\cdot)$ is
the Gamma function.
Then, for every $A,Z>0$,
$$
  \sup_{a\in S^{n-1}} \left|    \sum_{i=1} ^N(\langle X_i,  a\rangle^2   -
  \E\langle X_i,  a\rangle^2)\right|\leq {2A^2}{} + {6\sqrt n\,Z}{} +C_{\phi}
$$
with probability larger than
$$1- 4\exp(-n)-4\p(A_k>A)- 4\times 9^n\sup_{a\in S^{n-1}} \p \left(
\Big(\sum_{i>k} (\langle X_i,  a\rangle^*)^4 \Big)^{1/2}>Z\right).
$$
\end{lemma}

The term involving $Z$ in the upper bound will be bounded later using
general estimates in Lemma~\ref{lemma:remainder}.  Thus
Lemma~\ref{KLS technic} clearly stresses the fact that in order to
estimate the distance between the covariance matrix and the empirical
one, it will remain to estimate $A_k$, to get $A$.

\medskip

\noindent {\bf Proof:}
Let $\Lambda\subset \R^n$  be a (1/4)-net of the unit Euclidean ball in the Euclidean metric
of cardinality not greater than $9^n$. Let $(\varepsilon_i)_{1\leq i\leq N}$ be i.i.d. $\pm 1$
Bernoulli random variables of parameter $1/2$. By Hoeffding's inequality, for every $t>0$ and every
$(s_i)_{1\leq i\leq N}\in\R^N$,
$$
\p_{(\varepsilon_i)}\left( \Bigl| \sum_{i=1} ^N\varepsilon_is_i\Bigr|\geq
  t(\sum_{i=1}^N |s_i|^2)^{1/2}\right)\leq 2\exp(-t^2/2).
$$

Fix an arbitrary  $1\leq k<N$.  For every $(s_i)_{1\leq i\leq
  N}\in\R_+^N$  there exists a permutation
$\pi$ of $\{1,\ldots, N\}$  such that
$$
\Bigl| \sum_{i=1}^N \varepsilon_is_i \Bigr|\leq
\sum_{i=1}^k s_i^*+\Bigl| \sum_{i=k+1} ^N\varepsilon_{\pi(i)} s_i^*\Bigr|, 
$$
where 
$(s_i^*)_i$ denotes a non-increasing rearrangement of $(|s_i|)_i$.

Also,  it is easy to check using (\ref{eq:AkBk_intro})
that  for any $a\in S^{n-1}$ and any $I\subset \{1,\ldots,N\}$
with $|I|\leq k$, $\sum_{i\in I} \langle X_i, a\rangle^2\leq A_k^2$.

Thus, for  every $a\in S^{n-1}$,
$$
\p_{(\varepsilon_i)}\left( \Bigl| \sum_{i=1} ^N\varepsilon_i \langle X_i,
  a\rangle^2\Bigr|\leq A_k^2+t\Bigl(\sum_{i=k+1}^N (\langle X_i,
  a\rangle^*)^4\Bigr)^{1/2}\right)\geq 1- 2\exp(-t^2/2).
$$
Note that $\sum_{i=1} ^N\varepsilon_i \langle X_i,
  a\rangle^2= \sum_{i\in E}  \langle X_i,
  a\rangle^2
-\sum_{i\in E^c} ^N \langle X_i,
  a\rangle^2$ for some set $E\subset\{1, ..., N\}$ and we can 
apply a union bound argument indexed by $\Lambda$ together with 
Lemma~\ref{matrix}. We get that
$$
\p_{(\varepsilon_i)}\left(
\sup_{a\in S^{n-1}}\Bigl| \sum_{i=1} ^N\varepsilon_i
\langle X_i, a\rangle^2\Bigr|\leq 2
\Bigl[A_k^2+t
\sup_{a\in \Lambda} \Bigl(\sum_{i=k+1}^N (\langle X_i,
   a\rangle^*)^4\Bigr)^{1/2 }\Bigr]\right)
$$
$$
\geq 1- 2\times 9^n\exp(-t^2/2).
$$
Using again a union bound argument and the triangle inequality to
estimate the probability that the $(X_i)$ satisfy
$$
\sup_{a\in \Lambda} \Big(\sum_{i=k+1}^N (\langle X_i,
a\rangle^*)^4\Big)^{1/2}> Z,
$$
and choosing
$t=3\sqrt n$ (so that $2\cdot 9^n\exp(-t^2/2)\leq e^{-n}$)
 we get  that
$$
\sup_{a\in S^{n-1}}\Bigl| \sum_{i=1} ^N\varepsilon_i
\langle X_i,a\rangle^2\Bigr|\leq 2A^2+6\sqrt n\,Z
$$
with probability larger than
$$
1- e^{-n}-\p(A_k>A)- 9^n\sup_{a\in S^{n-1}} \p \left(
\Big(\sum_{i=k+1}^N (\langle X_i, a\rangle^*)^4 \Big)^{1/2}>Z\right).
$$

Now we transfer the result from Bernoulli random variables to centered
random variables (see \cite{LeTa}, Section~6.1). 
By the triangle inequality, for every $s,t>0$, one
has
$$
  m (s)\, \p\left(\sup_{a\in S^{n-1}}\Big|\sum_{i=1} ^N\big(\langle
    X_i, a\rangle^2-\E
  \langle  X_i, a\rangle^2 \big)\Big|>s+t\right)
$$
$$
   \leq 2\p\left(\sup_{a\in S^{n-1}}\Big| \sum_{i=1} ^N\varepsilon_i
   \langle X_i, a\rangle^2\Big|>t\right)
$$
where $m (s)= \inf_{a\in S^{n-1}} \p\left( \Big| \sum_{i=1}
  ^N\big(\langle X_i, a\rangle^2-\E \langle X_i, a\rangle^2 \big)\Big|
  \leq s\right)$.

To conclude the proof it is enough to find $s$ so that $m(s)\geq 1/2$.
To this end we will use a general Lemma~\ref{desymmetrization}
(below).  First consider $\phi(t)=t^p$. For $a\in S^{n-1}$, set
$Z_i=|\langle X_i, a\rangle|^2/\ta^{2/p}$ and $q=p/2$.  Then by
Lemma~\ref{desymmetrization} we have $m(s)\geq 1/2$ for $s= 4\ta^{2/p} N^{1/r}$
and $r=\min(p/2, 2)$.  Now consider $\phi (t)=(1/2)
\exp(t^{\alpha})$. Then for every $a\in S^{n-1}$ and every $i\leq N$
using hypothesis $\mbox{\bf H}(\phi)$ we have
$$
  \E |\langle X_i, a\rangle|^4 \leq 8 \ta \, \int _0^{\infty} t^3 \exp(-t^{\alpha})\, dt
  = \frac{8\ta }{\alpha} \, \Gamma\left(\frac{4}{\alpha}\r) :=\ta C_{\alpha} .
$$
Given $a\in S^{n-1}$, set $Z_i=|\langle X_i, a\rangle|^2/\sqrt{\ta C_{\alpha}}$.
Then $\E Z_i ^ 2 \leq 1$.  Applying again Lemma~\ref{desymmetrization} (with $q=2$),
we observe that $m(s)\geq 1/2$ for $s=4 \sqrt{C_{\alpha} N}$.  This completes the proof.
\qed

It remains to prove the following general lemma. For convenience
of the argument above, we formulate this lemma using two
powers $q$ and $r$ rather than just one.

\begin{lemma}
\label{desymmetrization}
Let  $q\geq 1$ and $Z_1,\dots, Z_N$ be  independent non-negative
random variables  satisfying
$$
   \forall 1\leq i\leq N\qquad \E Z_i^q\leq 1.
$$
Let $r=\min(q,2)$, then
$$\forall z\geq 4 N^{1/r}
\qquad
\p\left(  \Big|
\sum_{i=1} ^N\big(Z_i-\E Z_i \big)\Big|
\leq z\right)\geq \frac{1}{2}.$$
\end{lemma}

\noindent{\bf Proof:} By definition of $r$, we have for all $i =1,
\ldots, N,$ $\E Z_i^r \le 1$. Since the $Z_i$'s are independent, we
deduce by a  classical symmetrization argument that
\[
\E \Big|
\sum_{i=1} ^N\big(Z_i-\E Z_i \big)\Big|
\le 2 \E \E_{(\eps_i)} \Big|
\sum_{i=1} ^N\eps_i Z_i \Big|
\le 2 \E \left( \sum_{i=1} ^N Z_i^2 \right)^{1/2}\le 2 \E \left( \sum_{i=1} ^N Z_i^r \right)^{1/r}
\]
since $r \in [1,2]$.  From  $\E Z_i^r \le 1$, we get that
\[
\E \Big|
\sum_{i=1} ^N\big(Z_i-\E Z_i \big)\Big|
\le 2 \E \left( \sum_{i=1} ^N Z_i^r \right)^{1/r} \le 2  \left(
\sum_{i=1} ^N \E Z_i^r \right)^{1/r} \le  2 N^{1/r}.
\]
By Markov's inequality we get
\[
\p\left(  \Big|
\sum_{i=1} ^N\big(Z_i-\E Z_i \big)\Big|
\geq 4 N^{1/r} \right)\leq \frac{1}{2},
\]
and since $ z\geq 4 N^{1/r}$, this implies the required
estimate.
\qed

The following lemma is standard (cf. Lemma~5.8 in \cite{LeTa},
which however contains  a misprint).

\begin{lemma}
\label{lemma:remainder}
Let  $q>0$ and  let  $Z_1,\dots, Z_N$ be  independent non-negative
random variables  satisfying
$$
\forall 1\leq i\leq N\quad \forall t\geq 1\quad \p(Z_i\geq t)\leq
1/t^q.
$$
Then, for every $s>1$, with probability larger than $1-s^{-k}$, one
has
$$
\sum_{i=k}^N Z_i^*\leq
\begin{cases}
\frac{(2es)^{1/q}}{1-q}N^{1/q}k^{1-1/q} &\mbox{if } 0<q<1 \\
{2es}N\ln\Big(\frac{e N}{k}\Big) &\mbox{if } q=1 \\
\frac{12 q(es)^{1/q}}{q-1} N  & \mbox{if } q> 1.
\end{cases}
$$

\end{lemma}

\noindent{\bf Proof:} Assume first that $0<q\leq 1$.
It is clear that
$$\forall 1\leq i\leq N\qquad \p(Z_i^*>t)\leq \binom{N}{i}t^{-iq}\leq
(Ne/it^q)^i, $$
where  we used the inequality $\binom{N}{i}\leq (Ne/i)^i$.
Thus if $eNt^{-q}\leq 1$, then
$$\p(\sup_{i\geq k}i^{1/q}Z_i^*>t)\leq
\sum_{i\geq k} (Ne/t^q)^i= \Big(\frac{eN}{t^q}
\Big)^k(1- eNt^{-q})^{-1}.$$
Therefore if  $eNt^{-q}\leq 1/2$, then
$\p(\sup_{i\geq k}i^{1/q}Z_i^*>t)\leq
({2eN}{t^{-q}})^k.$
Since the inequality is trivially true if $eNt^{-q}\geq 1/2$,
it is proved for every $t>0$.
Therefore for $q<1$ we have
$$
  \sum_{i=k}^N Z_i^*\leq t\sum_{i=k}^{\infty} i^{-1/q}
  \leq t\left( k^{-1/q} - \frac{k^{1-1/q}}{1-1/q} \r) \leq
  \frac{t}{1-q} \, k^{1-1/q}
$$
with probability larger than $1-(2eN/t^q)^k$. Choosing
$t=(2esN)^{1/q}$, we obtain the estimate in the case $0<q<1$.

For $q=1$ we have
$$\sum_{i=k}^N Z_i^*\leq t\sum_{i=k}^{N} i^{-1}
\leq t\left( \frac{1}{k} + \ln (N/k) \r) \leq  t\,  \ln (e N/k)
  $$
with probability larger than $1-(2eN/t)^k$. To obtain the desire
estimate choose  $t=2esN$.

Now assume that $q>1$. Set $\ell=\lceil\log_2 k\rceil$.
The same computation as before for the scale $(2^{i/q})$ instead of
$(i^{1/q})$ gives that
$$
\p\left(\sup_{i\geq \ell}2^{i/q}Z_{2^i}^*>t\r)\leq
\sum_{i\geq \ell} \left(Net^{-q}\r)^{2^i}\leq\left({2eN}{t^{-q}}
\r)^{2^\ell}.
$$
Note also that
$$
  \p\left(k^{1/q} Z_{k}^*>t\r)\leq \left(N e t^{-q}\r)^{k}.
$$
Thus
$$
  \sum_{i=k}^N Z_{i}^*  \leq  k Z_k^* +
  \sum_{i=\ell}^{\lceil\log_2 N\rceil} 2^iZ_{2^i}^*
  \leq
    t \left(k^{1-1/q} +  (4N)^{1-1/q}/(2^{1-1/q}-1)\r)
$$
$$
   \leq t \left(k^{1-1/q} + \frac{2q}{1-q}  (4N)^{1-1/q}\r) \leq
    t \, \frac{3q}{1-q}  (4N)^{1-1/q}
$$
with probability larger than $(N e t^{-q})^{k} + (2N e t^{-q})^{k}$.
Thus, taking $t=(4esN)^{1/q}$, we obtain
$$\p\left(
\sum_{i=k}^N Z_{i}^*
\leq \frac{12 q (es)^{1/q}}{q-1} N \right)\geq 1-s^{-k}.
$$
\qed

We are now ready to tackle the problem of approximating the covariance matrix by
the empirical covariance matrices, under hypothesis $H(\phi)$ with $\phi(t)=t^p$.
As our proof works for all $p>4$, we also include the case $p >8$ originally solved
in \cite{MPaouris2012} (under additional assumption on $\max_i |X_i|$). For clarity, we
split the result into two theorems. The case $4 < p \leq 8$ has been stated as
Theorem~\ref{thm:KLS-M} in the Introduction.

Before we state our result, let us remark that $p>2$ is a necessary condition. Indeed,
let $(e_i)_{1\leq i\leq n}$ be an orthonormal  basis  of $\R^n$ and let $Z$ be a random vector
such that  $Z= \sqrt n e_i$  with probability $1/n$. The covariance matrix of $Z$ is the identity $I$.
Let $A$ be an $n\times N$  random matrix with independent columns distributed as $Z$. Note that if
$\|\frac{1}{N}AA^\top-I\|<1$ with some probability, then $AA^\top$ is invertible with the same
probability. It is known (coupon collector's problem) that $N\sim n\log n$ is needed to have
$\{Z_i\,:\, i\leq N\}=\{\sqrt ne_i\,:\, i\leq n\}$ with probability, say, $1/2$. Thus for vector
$Z$, the hypothesis $H(\phi)$, $\phi(t)=t^2$ is satisfied but $N\sim n\log n$ is needed for the
covariance matrix to be well approximated by the empirical covariance matrices with probability $1/2$.

We also would like to mention that we don't know how sharp the power $\gamma/p$ 
appearing in the bound below is. In particular, it is not clear if it can be improved to $1/2$.

\begin{theorem}
\label{thm:KLS}
Let $4<p\leq 8$ and $\phi (t)=t^p$.
Let  $X_1,\dots, X_N$ be independent random vectors in $\R^n$ satisfying
hypothesis ${\bf H}(\phi)$.  Let $\eps \leq \min\{1,
(p-4)/4\}$ and $\gamma = p-4-2\eps$.
Then with probability larger than
$$
   1- 8 e^{-n} - 2  \eps^{-p/2}\max\left\{N^{-3/2}, n^{-(p/4-1)}\r\}
$$
one has
$$
  \sup_{a\in S^{n-1}} \left| \frac{1}{N}   \sum_{i=1} ^N(\langle X_i,
  a\rangle^2 - \E\langle X_i,  a\rangle^2)\right|\leq C \left(\frac{1}{N}\,
  \max_{i\leq N} |X_i|^2 + C(p, \eps)\, \left(\frac{n}{N}\r)^{\gamma/p} \r) ,
$$
where
$$
   C(p, \eps) = (p-4)^{-1/2}\, \eps ^{-4(2+\eps)/p} .
$$
and $C$ is an absolute constant.
\end{theorem}

An immediate consequence of this theorem is the following corollary.

\begin{cor}
Under assumptions of Theorem~\ref{thm:KLS-new}, assuming additionally that
$\max _i |X_i|^2 \leq C n^{\gamma/p} N^{1-\gamma/p}$ with high probability,
we have with high probability
$$
  \sup_{a\in S^{n-1}} \left| \frac{1}{N}   \sum_{i=1} ^N(\langle X_i,
    a\rangle^2   -
   \E\langle X_i,  a\rangle^2)\right|\leq C_1\,  C(p, \eps)\,
 \left(\frac{n}{N}\right)^{\gamma /p},
$$
 where $C$ and $C_1$ are absolute positive constants.
\end{cor}

\begin{theorem}
\label{thm:KLS-new}
There exists a universal positive constant $C$ such that the  following holds.
Let  $p> 8$, $\alpha \in (0,2]$. Let $\phi$ and $C_{\phi}$ be either
$\phi (t)=t^p$ and $C_{\phi}=C$ or $\phi (t)=(1/2) \exp(t^{\alpha})$
and $C_{\phi}=\left(C/\alpha\r)^{2.5/\alpha}$.
Let  $X_1,\dots, X_N$ be independent random vectors in $\R^n$ satisfying
hypothesis ${\bf H}(\phi)$. 
In the case $\phi (t)=t^p$ we also define
$$
   p_0 = 8 e^{-n} + 2 \left( \frac{3p-8}{6(p-8)}\r)^{p/2}\,
   N^{-(p-8)/8}\, n^{-p/8}
$$
and in the case  $\phi (t)=(1/2) \exp(t^{\alpha})$, we assume
$N\geq (4/\alpha)^{8/\alpha}$ and define
$$
   p_0 = 8 e^{-n} +  \frac{1}{(10 N)^4}  \exp\left(\frac{4 n^{\alpha/2}}{
  (3.5\ln(2n))^{2\alpha}}\right) + \frac{N^2}{2\exp((2 n N)^{\alpha/4})} .
$$
Then in both cases with probability larger than $1-p_0$ one has
$$
   \sup_{a\in S^{n-1}} \left| \frac{1}{N}   \sum_{i=1} ^N(\langle X_i,  a\rangle^2
  -   \E\langle X_i,  a\rangle^2)\right|
    \leq  \frac{C}{N}\max _{i\leq N} |X_i|^2  + C_{\phi} \sqrt{\frac{n}{N}} .
$$
\end{theorem}

As our argument works in all cases we prove both theorems together.

\smallskip

\noindent
{\bf Proof of Theorems~\ref{thm:KLS} and \ref{thm:KLS-new}.}
We first consider the case $\phi = t^p$. Note that in this case
$$
   \E \, |\langle X_i, a\rangle|^4 \leq 1+ \int_1^{\infty} \p\left( |X_i|^4 >t\r)\, dt
   \leq 1 + \int_1^{\infty} 4 s^{3-p}\, ds = \frac{p}{p-4}.
$$
Thus, by  Lemma \ref{KLS technic} it is enough to estimate
$A^2 + \sqrt{n} \,Z + \sqrt{p/(p-4)}\, \sqrt{N}$ and the corresponding probabilities.
We choose $k=n$.

In the case $\phi = t^p$ we  apply Lemma~\ref{lemma:remainder} with
$Z_i=|\la X_i, a\ra|^4$, $i\leq N$,  $q=p/4>1$ and $s=9 e$.
It gives
$$
  \p \left(\left(\sum_{i>n} (\langle X_i,  a\rangle^*)^4 \r)^{1/2}  >
      Z \right) \leq (9 e)^{-n},
$$
for
$$
  Z= \sqrt{\frac{12q}{q-1}} (es)^{1/2q}\, \sqrt{N} =
  \sqrt{\frac{12p}{p-4}} (3 e)^{4/p}\, \sqrt{N} .
$$

Now we estimate $A_n$, using Theorem~\ref{A_k}.

\bigskip

\noindent
{\bf Case 1: $4< p\leq 8$ } (Theorem~\ref{thm:KLS}). \,
We apply Theorem~\ref{A_k} (and the Remark  following it),
with $\sigma = 2+\eps$, where
$\eps<  (p-4)/4$,  $\lam =3$ and $t= 3 N^{2/p} n^{\delta}$
for  $\delta = 1/2 - 2/p$.   Then
$$
   M_1  \leq C(p, \eps) \sqrt{n} \, (N/n)^{(2+\eps)/p} ,
$$
where
$$
  C_0 (p, \eps) = C \left(\frac{1}{p-4}\r)^{(p-4-2\eps)/p} \,
    \left(\frac{1}{\eps}\r)^{2(2+\eps)/p}\,
$$
and
$$
  \beta \leq \frac{1}{5}\, \left(\frac{12}{ 5 e \eps N}\r)^3  +
    \frac{1}{4 (p-4)^p n^{\delta p}} \leq
   \eps^{-p }\max\{N^{-3}, n^{-\delta p}\}\leq 1/64
$$
provided that $n$ is large enough.
Then, using $\delta = 1/2-2/p$, we obtain
$$
 A_n^2 \leq  C \left(\max _{i\leq N} |X_i|^2 + N^{2/p} n^{\delta} \max
   _{i\leq N} |X_i|  + C_0^2(p, \eps)\,   n \, (N/n)^{2(2+\eps)/p}  \r)
$$
$$
  \leq  2 C \left(\max _{i\leq N} |X_i|^2 + C_0^2(p, \eps) \,
  n \, (N/n)^{2(2+\eps)/p}  \r) .
$$
Combining all estimates and noticing that $(p-4)^{-\gamma}<2$,
we obtain that the desired estimate holds with probability
$$
   1- 8 e^{-n} - 2  \eps^{-p/2 }\max\{N^{-3/2}, n^{-(p/4-1)}\} .
$$

\bigskip

\noindent
{\bf  Case 2:  $p> 8$ } (Theorem~\ref{thm:KLS-new}). \,
In this case we apply Theorem~\ref{A_k}
(see also the Remark following it),
with $\sigma =p/4$, $\lambda =(p-4)/2$, $t=3 (nN)^{1/4}$.
Then $M_1 \leq C \sqrt{n} (N/n)^{1/4}$ and
$$
  \beta \leq \left( \frac{2(3p-8)}{5e (p-8) N}\r)^{(p-4)/2}
    \frac{1}{p-5} + \frac{(3p-8)^p}{4 (6(p-8))^p N^{(p-8)/4} n^{p/4}}
$$
$$
 \leq \left( \frac{3p-8}{6(p-8)}\r)^{p}\, N^{-(p-8)/4}\, n^{-p/4} \leq 1/64,
$$
provided that  $N$ is large enough. Thus with probability at least
$1-\sqrt\beta$ we have
$$
   A_n^2 \leq C \left(\max _{i\leq N} |X_i|^2 + (nN)^{1/4} \max
     _{i\leq N} |X_i|  + \sqrt{nN} \r)
  \leq 2 C \left( \max _{i\leq N} |X_i|^2 + \sqrt{nN}   \r).
$$
Combining all estimates we obtain that the desired estimate holds with
probability
$$
   1- 8 e^{-n} - 2 \left( \frac{3p-8}{6(p-8)}\r)^{p/2}\,
     N^{-(p-8)/8}\, n^{-p/8}  .
$$

\bigskip

\noindent
{\bf  Case 3: $\phi (t)=(1/2) \exp(t^{\alpha})$ } (Theorem~\ref{thm:KLS-new}). \,
As in Case 2 we apply  Lemma~\ref{KLS technic}. It implies that it is
enough to estimate $A^2 + \sqrt{n} \,Z + \sqrt{C(\alpha) N}$, with
$C(\alpha)$ from Lemma~\ref{KLS technic}, and the corresponding probabilities.
A direct calculations show that in this case we have for
$C_{\alpha}'=(4/\alpha)^{1/\alpha}$ and $t>1$,
$$
   \p\left( \left( |X|/C_{\alpha}' \r)^4 >t \r) \leq 2 \exp(C_{\alpha}') \,
  t^{\alpha/4}  \leq \frac{1}{t^2}.
$$
We  apply Lemma~\ref{lemma:remainder} with
$Z_i= |\la X_i, a\ra|^4 / \sqrt{C_{\alpha}'}$, $i\leq N$,  $q=2$ and
$s=9 e$.
It gives
$$
  \p \left(\left(\sum_{i>k} (\langle X_i,  a\rangle^*)^4 \r)^{1/2}  >
      Z \right) \leq (9 e)^{-n},
$$
for
$$
  Z= \left(C_{\alpha}'\r)^{1/4} 6 \sqrt{6} e\, \sqrt{N} .
$$

To estimate $A_n$ we use Theorem~\ref{A_k} with $t=(n N)^{1/4}$ and
$$
  \lam = 10 \left( N/n \r)^{\alpha/4} \min\left\{1, \left(\alpha \, \ln
        (2 N/n)\r)^{-1}\r\}.
$$
Note that
$$
    \max\left\{ 4,  10 \left( N/n \r)^{\alpha/4} \, \left(\ln
        (2 N/n)\r)^{-1}  \r\} \leq \lam \leq 10 \left( N/n \r)^{\alpha/4} .
$$
Then for absolute positive constants $C$, $C'$,
$$
  M_1 \leq  \sqrt{n}\,  \left(C \lam \r)^{1/\alpha}\,
   \left(  \ln \frac{2 N}{n}  +
     \frac{1}{\alpha} \r)^{1/\alpha}
  \leq  \left(\frac{C' }{\alpha}\r)^{1/\alpha}  \left(n N\r)^{1/4}
$$
and
$$
  \beta \leq  \frac{1}{(10 N)^4}  \exp\left(\frac{4 n^{\alpha/2}}{
   (3.5\ln(2n))^{2\alpha}}\right) + \frac{N^2}{2\exp((2 n N)^{\alpha/4})} \leq 1/64,
$$
provided that $N\geq (4/\alpha)^{8/\alpha}$. Thus with probability at least
$1- \sqrt\beta$ we have
$$
   A_n^2 \leq  C'' \, \max _{i\leq N} |X_i|^2 +
              \left(\frac{C'''}{\alpha}\r)^{2/\alpha} \sqrt{nN}  ,
$$
where $C''$ and $C'''$ are absolute positive constants.  This together
with the estimate for $Z$  completes the proof (note  that
$C(\alpha) \leq C(2/\alpha)^{5/\alpha}$).
\qed

\section{The proof of Theorem~\ref{A_k}}
\label{sec:proof_C}

In this section we prove the main technical  result of this paper,
Theorem~\ref{A_k}, which establishes upper bounds for norms of submatrices
of random matrices with independent columns.  Recall that
for  $1\leq k\leq N$ the parameters  $A_k$ and $B_k$ are defined by
(\ref{eq:AkBk_intro}).

\subsection{Bilinear forms of independent vectors}
\label{bilinear}

Let $X_1$,... $X_N$ be independent random vectors and $a\in
\R^N$. Given  disjoint
sets
$T, S\subset \{1, ..., N\}$ we let
\begin{equation}
  \label{qats}
  Q(a, T, S) = \left| \la \sum _{i\in T} a_i X_i, \sum _{j\in S} a_j X_j \ra \r|,
\end{equation}
with the convention that $\sum _{i\in \emptyset} a_i X_i=0$.

\medskip

The following two lemmas are in the spirit of Lemma~2.3 in \cite{MPaouris2012}.
Recall that  $(s_i^*)_i$ denotes a
non-increasing rearrangement of $(|s_i|)_i$.

\begin{lemma}
\label{split}
Let $X_1$,... $X_N$ be independent random vectors in $\R^n$. Let
$\gamma \in (1/2, 1)$,
 $I\subset \{1, ..., N\}$,  and $a\in \R^N$.
Let $k\ge|\supp (a)|$. Then there exists $\bar a\in \R^N$ such that
$\text{\rm supp}(\bar a)\subset \text{\rm supp}\, (a)$,
$|\supp (\bar a)|\leq \gamma k$,
$|\bar a|\leq |a|$,
and
$$
  Q(a, I, I^c)\leq Q(\bar a,I, I^c) +
  \max\left\{\sum _{ i=m }^{m+\ell -1} V_i^*, \sum _{ i=m }^{m+\ell -1} W_i^* \r\},
$$
where $\ell = \lceil (1-\gamma) k \rceil$, $m=\lceil (\gamma - 1/2) k
\rceil$,  and
$$
   V_i = \la  a_i X_i, \sum _{j\in I^c} a_j X_j \ra \quad \mbox{ for }  i\in I,
$$
$$
   W_j = \la \sum _{i\in I} a_i X_i,  a_j X_j \ra \quad \mbox{ for }  j\in I^c.
$$
\end{lemma}

\noindent
{\bf Proof.} Let $E\subset \{1, ..., N\}$ be such that $\supp
 (a) \subset E$ and $|E|=k$.  Everything is clear when $k=0$
or 1, because then $Q(a,I,I^c)=0$. Thus we may assume that $k\geq 2$.
Let $F_1=E\cap I$ and $F_2=E\cap I^c$.  First assume that
$s:=|F_1|\geq k/2$.  Note that $ (1-\gamma) k\leq k/2\leq s$, so that
$\ell\leq s$. 
Let
$J\subset F_1$   be a set with $|J|=\ell$ such that
the set $\{|V_j|: j\in J\}$ consists
of $\ell$ smallest values among the values  $\{|V_i|: i\in F_1\}$.
(That is,
$J\subset F_1$ is such that
$|J|=\ell$  and for all $j\in J$ and $i\in F_1\setminus J$
we have $|V_i|\ge |V_j|$.)
Now we let
$$
  \bar F_1= F_1\setminus J \quad \mbox{ and } \quad
  \bar F_2= F_2.
$$
Define  the vector $\bar a\in \R^N$ by the conditions
$$
\bar a_{|{\bar F_1}}=a_{|{\bar F_1}}, \qquad
\bar a_{|{\bar J}}= 0, \qquad
\bar a_{|{\bar F_2}}=a_{|{\bar F_2}}.
$$
Thus  $\bar a$ differs from
$a$ only on coordinates from $J$; in particular its support has
cardinality less than or equal to $|\supp (a)|-|J| = s-\ell \le k -\ell
= \gamma k$. Moreover,
\begin{align*}
  Q(a,I, I^c) &=
\left| \la \sum _{i\in {F_1} } a_i X_i, \sum _{j\in {F_2}} a_j X_j \ra \r|\\
&\leq \left| \la \sum _{i\in {J} } a_i X_i, \sum _{j\in {F_2}} a_j X_j \ra \r|
+ \left| \la \sum _{i\in {F_1\setminus J} } a_i X_i, \sum _{j\in {F_2}} a_j
X_j \ra \r|\\
&=   \left| \sum _{ i\in J }  \la  a_i X_i,
\sum _{j\in F_2} a_j X_j \ra \r|
+ Q(\bar a, I, I^c).
\end{align*}
Then we have
\begin{align*}
 Q(a,I, I^c)
&\leq Q(\bar a, I, I^c) +
  \sum _{ i\in J } \left| \la  a_i X_i,
\sum _{j\in F_2} a_j X_j \ra\right| \\
 &\leq
 Q(\bar a, I,I^c) + \sum _{i=s-\ell +1}^{s} V_i^*  \leq
  Q(\bar a, I,I^c) + \sum _{i=m}^{m+\ell - 1} V_i^*,
\end{align*}
where $m=\lceil (\gamma -1/2) k\rceil$ and using that $s-\ell+1\geq k/2 -
\lceil (1-\gamma) k\rceil+1> (\gamma -1/2) k$.

If $|F_1|< k/2$ then $|F_2|\geq k/2$ and we proceed similarly
interchanging the role of $F_1$ and $F_2$ and obtaining
$$
  Q(a, I,I^c)\leq
 Q(\bar a, I,I^c) +  \sum _{i=m}^{m+\ell - 1} W_i^* .
$$
\qed

\begin{lemma}
\label{star}
Let $\ta\geq 1$ and $X_1,\cdots,X_N$ be  independent random  vectors in $\R^n$ satisfying
hypothesis $\mbox{\bf H}(\phi)$ for some function $\phi\in \cal{M}$ with parameter $\ta$.
Let $a \in \R^N$ with $|a|=1$.
In the notation  of Lemma~\ref{split}, for every $t>0$ one has
$$
  \p \left(\sum _{i=m}^{m+\ell - 1} U_i^* > t A_k \r)\leq
     (2\ta)^k  \left(\phi \left( \frac{t \sqrt{m}}{\ell}
      \r) \r) ^{-m}
\le   (2\ta) ^k  \left(\phi \left( \frac{t \sqrt{\gamma _0 k}}{(1-\gamma) {k}+1}
 \r) \r) ^{-\gamma _0 k},
$$
where $\{U_i\}_i$ denotes either $\{V_i\}_i$ or $\{W_i\}_i$,
     and $\gamma_0=\gamma -1/2$.
\end{lemma}

\noindent
{\bf Remarks. 1. }
Taking $\phi (t) = t^p$ for some $p>0$, we obtain that if
\begin{equation}\label{cond4mom}
    \p\left( |\la X_i, a\ra |\geq t\r) \leq t^{-p}
\end{equation}
then
\begin{equation}\label{est4mom}
  \p \left(\sum _{i=m}^{m+\ell - 1} U_i^* > t A_k \r)\leq
    (2\ta)^k \left( \frac{t\sqrt{m}}{\ell}
    \r) ^{-mp}.
\end{equation}
Note that the condition (\ref{cond4mom}) is satisfied if
$$
   \sup _{i\leq N} \sup _{a\in S^{n-1}}\E |\la X_i, a\ra |^p \leq \ta.
$$
{\bf 2. } Taking $\phi = (1/2)\exp(x^{\alpha})$ for some $\alpha>0$, we obtain that if
\begin{equation}\label{condpsialpha}
     \p\left( |\la X_i, a\ra |\geq t\r) \leq 2\exp(-t^{\alpha})
\end{equation}
then
\begin{equation}\label{estpsialpha}
  \p \left(\sum _{i=m}^{m+\ell - 1} U_i^* > t A_k \r)\leq
    (2\ta)^{k+m} \exp\left( -m \left(\frac{t \sqrt{m}}{\ell}
    \r) ^{\alpha}\r).
\end{equation}
Note that the condition (\ref{condpsialpha})  is satisfied if
$$
   \sup _{i\leq N}  \sup _{a\in S^{n-1}}\E
   \exp\left(|\la X_i,
   a\ra | ^{\alpha}\r) \leq 2\ta.
$$

\noindent
{\bf Proof.}  Without loss of generality assume that $U_i=V_i$ for every $i$. 
Then
$$
 \sum _{i=m}^{m+\ell - 1} V_i^*  \leq  \ell V_m^* .
$$
Let $F_1 = \supp (a)\cap I $ and $F_2 = \supp (a)\cap I^c $.
Note that $V_m^*>s$ means that there exists a set $F\subset F_1$
of cardinality $m$ such that $V_i > s$ for every $i\in F$  (if cardinality
of $F_1$ is smaller than $m$, the estimate for probability is trivial).
Since $|F_1|\leq k$, we obtain
$$
  \p \left(\sum _{i=m}^{m+\ell-1} V_i^* > t A_k \r)\leq
  \p \left(\ell V_m^* > t A_k \r) \leq {k \choose m}
  \max _{F\subset F_1 \atop |F|=m}
  \p \left(\forall i\in F:\, \,  |V_i| > \frac{t A_k}{\ell} \r).
$$
Denote $Z:=\sum _{j\in F_2} a_j X_j$.
Since $|a|\le 1$ then
$|Z|\leq A_k$, and
note
that the $X_i$'s, $i\in F_1$ are independent of
$Z$. Thus, conditioning on $Z$ we obtain
\begin{align*}
  \p \left(\sum  _{i=m}^{m+\ell-1} V_i^* > t A_k\r)
    &  \leq 2^k  \max _{F\subset F_1 \atop
  |F|=m}\prod _{i\in F}\p \left( |a_i| |\la X_i, Z\ra |
     >  \frac{t A_k}{\ell}  \r) \\
   &   \leq (2\ta)^k  \max _{F\subset F_1 \atop |F|=m}\prod _{i\in F}
        \left(\phi \left(\frac{t}{\ell |a_i|}\r)\right)^{-1}.
\end{align*}

Now we show that for every $s>0$,
$$
  \prod _{i\in F} \left(\phi \left(\frac{s}{|a_i|}\r)\right)^{-1}
 \leq \left(\phi \left(
      s\sqrt{m} \r)\r)^{-m} ,
$$
Indeed, this estimate  is equivalent to
$$
  \frac{1}{m}\, \sum_{i\in F} \ln \phi \left(\frac{s}{|a_i|}\r) \geq 
   \ln \phi \left( s \sqrt{m} \r) ,
$$
which holds by convexity of $\ln \phi (1/\sqrt{x})$, the
facts that $|a|\leq 1$ and $|F|=m$, and since $\phi$ is
increasing. Taking $s=t/\ell$, we obtain
$$
  \p \left(\sum  _{i=m}^{m+\ell-1} V_i^* > tA_k \r) \leq (2\ta)^k
  \left(\phi \left( t\, \sqrt{m}/\ell \r)\r) ^{-m}.
$$
Finally note that $m=\lc (\gamma -1/2) k\rc \geq \gamma _0 k$ and
$\ell =\lc (1-\gamma) k\rc \leq (1-\gamma) k+1$.  Since $\phi$ is
increasing, we obtain the last inequality,  completing  the proof.
\qed

\subsection{Estimates for  off-diagonal of bilinear forms}
\label{bilin_indep}

For $1\le k\leq N$ and $I\subset \{1, ..., N\}$ we  define $Q_k(I) $ by
\begin{equation}
\label{qki}
   Q_k(I) = \sup _{E\subset \{1, ..., N\}\atop |E|\leq k} \sup _{a\in B_2^E}
   Q(a, E\cap I, E\cap I^c).
\end{equation}

Lemmas \ref{split}, \ref{star} and \ref{matrix} imply the following
proposition.

\begin{prop} \label{iter} Let $\ta \geq 1$ and 
  $X_1,\cdots,X_N$ be  independent random vectors in $\R^n$
  satisfying hypothesis $\mbox{\bf H}(\phi)$ with parameter $\ta$ for some function
  $\phi \in \cal{M}$.  Let $\eps\in (0, 1/2)$, $2\le k \le N$,
  $I\subset \{1, ..., N\}$, $\gamma \in (1/2, 1)$, and $\gamma_0 =
  \gamma -1/2$.  Then for every $t>0$ one has
$$
  \p \left( Q_k(I) > \frac{Q_{[\gamma k]}(I) + t A_k}{1-2\eps} \r)\leq
  \exp\left(k\left(\ln \frac{5\ta  e  N}{k\eps} -
  \gamma_0
   \ln\phi \left( \frac{t \sqrt{\gamma_0 k}}{(1-\gamma) {k}+1}
    \r)  \r)\r).
$$
Moreover, letting  $M=\max _i |X_i|$  one has, for all $\ell >1$ and $t>0$,
$$
  \p \left(Q_{\ell} (I) > t M \r)\leq  \frac{N^2\ta}{4\phi(4t/\ell)} .
$$
\end{prop}

\noindent
{\bf Proof. } For every $E\subset {1,..., N}$ with $|E|=k$
let ${\cal{N}}_E$ be an $\eps$-net in $B_2^E$
of cardinality at most $\left(2.5/\eps\r)^k$. 
Let $\cal{N}$ denote the union of ${\cal{N}}_E$'s.
Lemma~\ref{matrix} yields
$$
   Q_k(I) \leq (1-2\eps) ^{-1}
    \sup _{E\subset \{1, ..., N\}\atop |E|\leq k} \sup _{a\in {\cal{N}}_E}
   Q(a, E\cap I, E\cap I^c).
$$
Therefore, applying Lemmas~\ref{split} and \ref{star}, we observe that the event
$$
   Q_k(I) \leq (1-2\eps) ^{-1} \left(
    \sup _{E\subset \{1, ..., N\}\atop |E|\leq \gamma k} \sup _{a\in {\cal{N}}}
   Q(a, E\cap I, E\cap I^c) + t A_k\r)
$$
occurs with probability at least
$$
  1 - {N\choose k} \, \left(\frac{2.5}{\eps}\r)^k\,
  (2\ta)^k
  \left(\phi \left( \frac{t \sqrt{\gamma _0 k}}{(1-\gamma) {k}+1}
    \r) \r) ^{-\gamma_0 k}
   $$
This implies the first estimate.

Now we prove the ``moreover" part. For every $E\subset \{1, \ldots,
N\}$ of cardinality $\ell$ denote $F_1=E\cap I$, $F_2=E\cap I^c$,
$m=|F_1|$ (so $|F_2|=\ell -m$). We also denote
$$
    M_0 :=\max_{i\in I} \max _{j \in I^c} |\la X_i, X_j\ra | \quad
    \mbox{ and } \quad
     M_1 := \max _{j \in I^c} | X_j |
$$
Then for any $a\in B_2^E$ we have
$$
  \left|\la \sum _{i\in F_1} a_i X_i, \sum _{j \in F_2} a_j X_j\ra \r|  \leq
   \left|\sum _{i\in F_1} a_i\,  \sum _{j \in F_2} a_j \r| \, M_0
$$
$$
 \leq
  \sqrt{m (\ell -m )} \left(\sum _{i\in F_1} a_i^2 \r) ^{1/2}\, \left(\sum _{j\in F_2} a_j^2 \r) ^{1/2} M_0
  \leq \frac{\ell }{2} \, \frac{M_0}{2}.
$$
Therefore, by  the  union bound,
\begin{align*}
    \p \left(Q_{\ell} (I) > t M_1 \r)& \leq \p \left(M_0 > 4 t M_1
        /\ell \r) \\
& \leq
    \sum _{i\in I} \sum _{j \in I^c} \p \left( |\la X_i, X_j\ra | > 4
      t M_1 /\ell \r).
\end{align*}
Finally, using the fact that $X_i$ is independent of $X_j$ for $i\ne
j$, $|X_j|\leq M_1$ for every $j\in I^c$,
 and using  the tail behavior of variables $\la X_i, z\ra$, we obtain
$$
    \p \left(Q_{\ell} (I) > t M \r) \leq \p \left(Q_{\ell} (I) > t M_1
        \r)\leq
   \frac{|I|\, |I^c|}{\phi (4t/\ell)}  \leq  \frac{N^2\ta }{4 \phi (4t/\ell)} .
$$
\qed

\medskip

\begin{prop}
\label{estforq}
Let $1\leq k\leq N$.  Let $\ta \geq 1$ and 
  $X_1,\cdots,X_N$ be  independent random vectors in $\R^n$
  satisfying hypothesis $\mbox{\bf H}(\phi)$ with parameter $\ta$ 
 for some function $\phi \in \cal{M}$.
Let $t>0$, $\lambda\geq 1$.

\smallskip

\noindent{\bf Case 1. } Let $p>4$ and $\phi(x)=x^p$. Let $\sigma \in (2, p/2)$.
Then 
$$
   Q_k(I)\leq
   e^4 \, \left( t \max_{i\leq N} |X_i| +
  C_2(\sigma,\lambda,p)\sqrt{k}
  \left(\frac{5\ta eN}{k}\right)^{\sigma/p}
  \,  A_{k} \r)
$$
occurs with probability at least
\begin{equation}
\label{probest3}
    1-  \left(\frac{2(\sigma+\lambda)}{5\ta eN(\sigma-2)}\right)^{\lambda}
\frac{1}{2\lambda-1}  - \frac{N^2 \ta (\sigma +\lambda)^p}{4  (2t
  (\sigma -2))^p}
\end{equation}
and
$$
  C_2(\sigma,\lambda,p)=8\,
  \sqrt{\frac{\sigma+\lambda}{1+\lambda/2}}\,
  \left(\frac{2p}{p-2\sigma}\right)^{1+2\sigma /p}
  \left(\frac{2(\sigma+\lambda)}{\sigma-2}\right)^{2\sigma/p}.
$$

\noindent{\bf Case 2. } Assume that $\phi(x)=(1/2) \exp(x^\alpha)$ for some
$\alpha>0$.
Then for every $t>0$,
$$
 Q_k(I)\leq C^{1/ \alpha}  \left( t \max_{i\leq N} |X_i| + \left(C \lambda \r)^{1/\alpha}
  \sqrt{k} \left( \left(\ln \frac{20 \ta e N }{k}\r)^{1/ \alpha} +
  \left(\frac{1}{\alpha}\r)^{1/ \alpha} \r)  A_k\r)
$$
 with probability at least
$$
  1- \frac{1}{(10 \ta N)^{\lambda}}  \exp\left(- \frac{\lambda k^{\alpha/2}}{(3.5\,
  \ln(2 k))^{2\alpha}}\r)  -  \frac{N^2 \ta}{2\exp((2t)^{\alpha})}.
$$
\end{prop}

\noindent
{\bf Proof.} Let $\gamma\in (1/2,1)$ to be chosen later.
 For  integers $s\geq 0$ denote $k_0=k$, $k_{s+1}= [\gamma  k_s]$.
Clearly, the sequence is strictly decreasing whenever $k_s\geq 1$
and $k_s\leq \gamma^s k$. Assume that $k\geq 1/(1-\gamma)$.
Define $m$ to be the largest integer $m\geq 1$ such that
$ k_{m-1}\geq 1/(1-\gamma).$ Note that $\gamma k_{m-1}\geq 1$.
Therefore
\begin{equation}\label{defkm}
  1\leq k_m< \frac{1}{1-\gamma} \leq k_{m-1}.
\end{equation}

By Proposition~\ref{iter} we observe that for every positive $t_{s}$
and $\varepsilon_s\in (0,1/2)$,
$0\leq s\leq m$, the event
$$
   Q_k(I)\leq
   \left( Q_{k_m}(I) + \sum _{s =0}^{m-1} t_{s} A_{k_{s}} \r)\
   \prod _{s =0}^{m-1} (1-2\eps _{s})^{-1}
$$
occurs with probability at least
\begin{equation}
\label{probest1}
   1- 2\sum _{s =0}^{m-1}\exp\left(k_{s}
   \left(\ln \frac{5\ta  e N}{k_{s} \eps_{s} } -
   \gamma_0 \,
   \ln\phi \left( \frac{t_s \sqrt{\gamma_0 k}}{(1-\gamma) {k}+1}
    \r)   \r) \r).
\end{equation}
Let $\eps >0$ and a positive decreasing sequence $(\eps _s)_s$
be chosen later and set
$$
   t_s= \frac{(1-\gamma)k_s+1}{\sqrt{\gamma _0 k_s}}\
   \phi^{-1}\left(\left(\frac{5\ta eN}{ k_s\varepsilon_s}
   \right)^{(1+\varepsilon)/\gamma_0}\right),
$$
where $\phi^{-1}(s)=\min\{t\geq 0\,:\, \phi(t)\geq s\}$.

We start estimating $Q_k(I)$. Since $\ln (1-x)\geq -2x$ on $(0, 3/4]$,
we observe that for $\eps _s<3/8$,
$$
 \sum _{s =0}^{m-1}\ln  (1-2\eps _{s}) \geq  \sum  _{s =0}^{m-1}
 - 4 \eps _{s}
$$
so that
$$
  \prod _{s =0}^{m-1} (1-2\eps _{s})^{-1}\le
  \exp\left( 4 \sum  _{s =0}^{m-1} \eps _{s} \r)
$$
Note that
$$
 \sum_{s=0}^{m-1} t_s A_{k_s} \leq A_k\sum_{s=0}^{m-1}t_s .
$$
Thus by (\ref{probest1}) and by our choice of $t_s$,
\begin{equation}\label{estforQ_k-1}
   Q_k(I)\leq \exp\left( 4 \sum  _{s =0}^{m-1} \eps _{s} \r)\,
   \left( Q_{k_m}(I) + A_k \sum _{s =0}^{m-1} t_s\r)\,
\end{equation}
with probability at least
$$
   1-  2\sum _{s =0}^{m-1}\exp\left(- k_{s}\, \varepsilon\, \ln
   \frac{5\ta e N}{k_{s}\varepsilon_s}   \r) \geq  1- 2  \exp\left(- k_{m-1}\,
   \varepsilon\, \ln \frac{5\ta e N}{k_{m-1}}  \r) \sum_{s=0}^{m-1} \eps_s^{k_s \varepsilon}.
$$
 Since $k_{m-1} \geq 1/(1-\gamma)$, this probability is larger than
\begin{equation}
\label{probest1-1}
    1- 2  \exp\left(- \frac{\varepsilon}{1-\gamma}\, \ln \left(5\ta e (1-\gamma) N \r) \r)
  \sum_{s=0}^{m-1} \eps_s^{k_s \varepsilon}.
\end{equation}
Thus it is enough to choose appropriately $\eps _s$ and to estimate
$\sum_{s=0}^{m-1}t_s$, $Q_{k_m}(I)$ and $\sum_{s=0}^{m-1} \eps_s^{k_s \varepsilon}$.
We distinguish two cases for $\phi$.

\medskip
\noindent {\bf Case 1: $\phi(x)=x^p$. \quad }
In this case we choose $\eps _s = (s+2)^{-2}$ so that
$$
   \sum_{s=0}^{m-1} \eps_s^{k_s \varepsilon} =  \sum_{s=0}^{m-1}
   (s+2)^{-2 k_s\varepsilon}\leq \sum_{s=0}^{m-1}(s+2)^{-2 k_{m-1}
    \varepsilon}\leq \frac{1}{2 k_{m-1}\varepsilon-1} .
$$
 Choose $\varepsilon=\lambda(1-\gamma)$. Since $\lam \geq 1$ and
$k_{m-1} \geq 1/(1-\gamma)$, we have
$2k_{m-1}\varepsilon\geq 2\varepsilon/(1-\gamma)= 2\lambda$
and
$$
  \sum_{s=0}^{m-1}(s+2)^{-2k_s\varepsilon}\leq \frac{1}{2\lambda-1}.
$$
Using again $k_{m-1}\geq (1-\gamma)^{-1}$, we conclude that the probability in
(\ref{probest1-1}) is larger  than
\begin{equation}
\label{probest2}
1- \left(5\ta e N(1-\gamma)\right)^{-\lambda}   \frac{2}{2\lambda-1}.
\end{equation}

Now we estimate $\sum_{s=0}^{m-1}t_s$. We have
$$
  t_s= \frac{(1-\gamma)k_s+1}{\sqrt{\gamma _0 k_s}}
    \phi^{-1}\left(\left(\frac{5\ta eN}{ k_s\varepsilon_s}
   \right)^{(1+\varepsilon)/\gamma _0} \right) =
   \frac{(1-\gamma)k_s+1}{\sqrt{\gamma_0 k_s}}
    \left(\frac{5\ta eN}{ k_s\varepsilon_s}
   \right)^{(1+\varepsilon)/\gamma_0 p}   .
$$
Recall that $\gamma> 1/2$, $k_{m-1} \geq 1/(1-\gamma)$,
so that  $(1-\gamma) k_s+1\leq 2(1-\gamma) k_s$
for $s\leq m-1$. Thus
$$
  t_s\leq \frac{2(1-\gamma)\sqrt k_s}{\sqrt{\gamma _0 }}
  \left(\frac{5\ta eN}{ k_s\varepsilon_s}\right)^{(1+\varepsilon)/\gamma _0 p}  .
$$
Let $b=(1+\varepsilon)/\gamma_0 p$. Assume that $b<1/2$. Since $k_s\leq \gamma ^s k$, we have
\begin{equation}\label{sum t_s}
  \sum _{s =0}^{m-1} t_{s} \leq
  \frac{2 (1-\gamma)  k^{1/2-b}(5\ta eN)^{b}}{\sqrt{\gamma_0}}\,
  \sum _{s =0}^{m-1} (s+2)^{\delta b}
  \gamma^{s(1/2-b)}.
\end{equation}
 Since  the function $h(z)=z^{2 b}{\gamma^{z(1/2-b)}}$ on $\R^+$ is
first increasing and then decreasing, we get
$$
 \sum _{s =0}^{m-1} (s+2)^{2 b} \gamma^{s(1/2-b)}
  =
  \gamma ^{-2(1/2-b)} \sum _{s =2}^{m+1}h(s) \leq
   \gamma ^{-1} \left(\sup_{z>0} h(z)  +  \int_0^\infty h(z)\, dz\r)
$$
$$
  \leq 2  \left(  \left( \frac{2 b}{(1/2-b) e \ln (1/\gamma)} \r)^{2 b}
   +
 \frac{\Gamma(1+2 b)}{((1/2-b)\ln(1/\gamma))^{1+2 b}} \r).
$$
As $2 b \leq 1$,  $\Gamma(1+2 b)\leq 1$.
Using also that $ \ln(1/\gamma)\geq 1-\gamma$, we observe that
 the previous quantity does not exceed
$$
   \frac{4}{((1/2-b)(1-\gamma))^{1+2 b}}.
$$
Coming back to (\ref{sum t_s}), we get
\begin{equation}\label{sumts}
 \sum _{s =0}^{m-1} t_{s} \leq
  \frac{ 8  k^{1/2-b} (5\ta eN)^{b}}{(1/2-b)^{1+2 b}
 (1-\gamma)^{2 b}\, \sqrt{\gamma -1/2}}  .
\end{equation}
To conclude this computation, we choose the parameter
$$
  \gamma=\frac{1+\lambda+\sigma/2}{\sigma +\lambda}.
$$
Note that $\gamma \in (1/2, 1)$ as required,
since $\lambda \geq 1$ and $2<\sigma$. With such a choice of $\gamma$,
we have
 $b=\sigma/p <1/2$, since $\sigma<p/2$.
Thus from (\ref{sumts}) and  (\ref{probest2})
$$
  \sum _{s =0}^{m-1} t_{s} \leq 8 \sqrt{k} \left(\frac{5\ta eN}{k}\right)^{\sigma/p}
  \left(\frac{p}{p/2-\sigma}\right)^{1+2\sigma /p}
  \left(\frac{\sigma +\lambda}{\sigma/2-1}\right)^{2\sigma/p}
  \sqrt{\frac{\sigma +\lambda}{1+\lambda /2}}
$$
holds with probability larger than
$$
  1- \left(5\ta e N\,\frac{\sigma/2-1}{\sigma+\lambda}\right)^{-\lambda}
  \frac{2}{2\lambda-1}.
$$

Finally, to estimate $Q_{k_m}$,  we note that
$$
   k_m < \frac{1}{1-\gamma}=\frac{\sigma+\lambda}{\sigma/2-1},
$$
and apply ``moreover" part of Proposition~\ref{iter}
(with $\ell = k_m$). Note that at the beginning of the proof
we assumed that $k\geq 1/(1-\gamma)$.
In the case   $k<1/(1-\gamma)$ the result  trivially holds
by the ``moreover" part of Proposition~\ref{iter} applied  with  $\ell=k$.

\medskip
\noindent {\bf Case 2: $\phi(x)=(1/2)\exp(x^\alpha)$.\quad }
   In this case we choose $\gamma=2/3$, so that $\gamma _0=1/6$.
As before we assume that $k\geq 1/(1-\gamma) = 3$ (otherwise $Q_k(I)\leq Q_2(I)$).
By (\ref{defkm}) we have $k_m < 3$,
hence, by (\ref{estforQ_k-1})
$$
   Q_{k}(I) \leq  \exp\left( 4 \sum  _{s =0}^{m-1} \eps _{s} \r)\,
   \left( Q_2(I) +  A_k \, \sum_{s=0}^{m-1} t_s  \r) .
$$

We define $k_s$ by 
$$
  \eps _s =\frac{1}{2}\,  \exp\left(-\left(\frac{k}{k_s}\r)^{\alpha/2}\, 
   \frac{1}{(s+2)^{2\alpha}} \r).
$$
Observe that since $k_s \leq \gamma ^s k$ and $\gamma = 2/3$, one has
$$
    \eps _s \leq \frac{1}{2}\, \exp\left(-\left(\frac{3}{2}\r)^{\alpha s/2}\,
    \frac{1}{(s+2)^{2\alpha}} \r) \leq \frac{1}{2e}\,  (s+2)^{2\alpha}
   \, \left(\frac{2}{3}\r)^{s \alpha/2} ,
$$
which implies
\begin{equation}\label{estfore_s-1}
   \sum_{s=0}^{m-1} \eps _s \leq \frac{C}{\alpha},
\end{equation}
for a positive absolute constant $C$.

We have
$$
   t_s=  \sqrt{6}\, \frac{k_s/3+1}{\sqrt{ k_s}}
    \phi^{-1}\left(\left(\frac{5\ta eN}{ k_s\varepsilon_s}
   \right)^{6(1+\varepsilon)} \right) =
  \sqrt{6}\,  \frac{k_s/3+1}{\sqrt{ k_s}}
    \left(\ln \left(2\left(\frac{5\ta eN}{ k_s\varepsilon_s}\r)^{6(1+\varepsilon)}\right)
   \right)^{1/ \alpha}   .
$$
By (\ref{defkm}) we have $k_m < 3\leq k_{m-1}$, hence,
\begin{align*}
   t_s  &\leq \sqrt{6}\, \frac{2}{3}\,
   \sqrt{ k_s}\, \left(6\left(1+\varepsilon\r)\r)^{1/\alpha}\,
   \left( \ln \frac{20\ta e N}{ k_s\varepsilon_s} + \ln\frac{1}{2 \varepsilon_s}
  \right)^{1/ \alpha}  \\
    &\leq  \sqrt{6}\, \frac{2}{3}\, 2^{1/\alpha}\, \sqrt{ k_s}\,
    \left(6\left(1+\varepsilon\r)\r)^{1/\alpha}\,
  \left(\left( \ln \frac{20\ta e N}{ k_s\varepsilon_s}\right)^{1/ \alpha}  + \left(
  \ln\frac{1}{2 \varepsilon_s} \right)^{1/ \alpha} \right).
\end{align*}
By the  choice of $\eps _s$ we obtain
\begin{equation}
\label{newcase2}
  \sum_{s=0}^{m-1}  \sqrt{ k_s}\, \left(\ln \frac{1}{2 \varepsilon_s} \right)^{1/ \alpha} \leq
  \sqrt{k} \, \sum_{s=0}^{m-1}  (s+2)^{-2} \leq 3 \sqrt{k}.
\end{equation}
Since $3^{-s}k \leq k_s\leq (2/3)^s k$, we observe
$$
  \sum_{s=0}^{m-1}  \sqrt{ k_s}\, \left( \ln \frac{20\ta e N}{ k_s\varepsilon_s}
   \right)^{1/ \alpha}  \leq \sqrt{ k}\, \sum_{s=0}^{m-1}  \left(\frac{2}{3}\r)^{s/2}
   \, \left(\ln \frac{20\ta e N 3^s}{k}\r)^{1/ \alpha}
$$
$$
  \leq \sqrt{k}\left( \sum_{s=0}^{m-1}  \left(\frac{2}{3}\r)^{s/2}
   \, 2^{1/\alpha}\, \left(\ln \frac{20 \ta e N }{k}\r)^{1/ \alpha} +
   \sum_{s=0}^{m-1}  \left(\frac{2}{3}\r)^{s/2}
   \, (2 s \ln 3)^{1/\alpha}\, \r)
$$
$$
 \leq  C_1^{1/\alpha} \sqrt{k} \left( \left(\ln \frac{20\ta e N }{k}\r)^{1/ \alpha} +
  \Gamma (1+1/\alpha) \r),
$$
where $C_1$ is an absolute positive constant and $\Gamma $ is the Gamma function.
This together with  (\ref{newcase2}) implies that
\begin{equation}\label{estfort_s-1}
  \sum_{s=0}^{m-1}  t_s \leq  \left(C_2(1+\eps)\r)^{1/\alpha} \sqrt{k}
  \left( \left(\ln \frac{20 \ta e N }{k}\r)^{1/ \alpha} +
  \Gamma (1+1/\alpha) \r),
\end{equation}
where $C_2$ is an absolute positive constant.

Now we estimate the probability. By the choice of $k_s$ we have
\begin{align*}
  \sum_{s=0}^{m-1} \eps_s^{\eps k_s}  &=  \sum_{s=0}^{m-1} \exp\left(-\eps k_s \ln (1/\eps _s)\r)
  = \sum_{s=0}^{m-1} \exp\left(-\eps k_s \left(\ln 2 + (k/k _s)^{\alpha/2} (s+2)^{-2\alpha}\r)\r) \\
  &\leq \sum_{s=0}^{m-1} \exp\left(-\eps k_s^{1-\alpha/2}\, k^{\alpha/2} (s+2)^{-2\alpha}\r).
\end{align*}

Since $k_s\geq k_{m-1}\geq 1/(1-\gamma)$ and $s+2\leq m+1$ for every $s\leq m-1$,
we get that
$$
  \sum_{s=0}^{m-1} \eps_s^{\eps k_s}\leq m \, \exp\left(-\frac{\eps}{(1-\gamma)^{1-\alpha/2}}
  \, \frac{k^{\alpha/2}}{(m+1)^{2\alpha}}\r) .
$$
Since $m$ is chosen such that $ 1/(1-\gamma)\leq k_{m-1}\leq  (2/3)^{m-1} k$, we observe
that
$$
   m-1 \leq \frac{\ln(k(1-\gamma))}{\ln(3/2)}.
$$
Therefore,
\begin{align*}
  \sum_{s=0}^{m-1} \eps_s^{\eps k_s}   &\leq \left(1+\frac{\ln(k/3)}{\ln(3/2)} \r)
  \exp\left(- \, \frac{\eps}{(1/3)^{1-\alpha/2}}\, \frac{k^{\alpha/2}}{(2.5\, \ln k)^{2\alpha}}\r)\\
  &\leq 2 \exp\left(- 3 \eps \frac{k^{\alpha/2}}{3^{\alpha/2}\, (2.5\, \ln k)^{2\alpha}}\r),
\end{align*}
which shows that probability in (\ref{probest1-1}) is at least
$$
    1- \frac{4}{(15\ta e N)^{3\eps}}
 \exp\left(- 3 \eps \frac{k^{\alpha/2}}{(3.5\, \ln k)^{2\alpha}}\r).
$$

Finally, to estimate $Q_2(I)$ we apply the ``moreover" part of Proposition~\ref{iter}
(with $\ell =2$). Choosing $\eps = \lambda/3$ and combining estimates (\ref{estfore_s-1}),
and (\ref{estfort_s-1}) with the estimate for $Q_2(I)$ we obtain the desired result.
\qed

\subsection{Estimating $A_k$ and $B_k$}
\label{subsec:A_k}

We are now ready to pass to the  proof of Theorem~\ref{A_k}.
To prove the theorem we  need two simple lemmas.

\begin{lemma}
\label{fubini} Let $\beta \in (0, 1)$.
Let $\p_1$ and $\p_2$ be probability measures  on $\Omega_1$ and $\Omega_2$
respectively and let $V\subset \Omega_1\otimes\Omega_2$ be such that
$$
  \p_1\otimes\p_2(V)\geq 1-\beta.
$$
Then there exists $W\subset \Omega_2$ such that
$$
  \p_2(W)\geq 1-\sqrt{\beta} \quad \mbox{ and }\quad
   \forall x_2 \in W,   \, \, \p_1\left(\left\{x_1:\, (x_1,x_2)\in V\r\} \r)\geq 1-\sqrt{\beta} .
$$
\end{lemma}

\noindent
{\bf Proof. }
Fix some $\delta \in (0,1)$.
Let
$$W:=\{x_2\in\Omega_2\,:\, \p_1\left(\left\{x_1\in\Omega_1:\, \,  (x_1,x_2)\in V\r\}\r)  \geq 1-\delta\}.$$
Clearly,
$$W^c=\{x_2\in\Omega_2\,:\, \p_1\left(\left\{x_1\in\Omega_1:\, \,  (x_1,x_2)\in V^c\r\}\r)  \geq \delta\}.$$
Then
$$
  \beta \geq \p_1\otimes\p_2(V^c) = \int_{\Omega_2}
  \p_1\left(\left\{x_1\in\Omega_1:\, \,  (x_1, x_2)\in V^c \r\}\r)\, d\,\p_2(x_2)
$$
$$
  \geq \int_{W^c} \p_1\left(\left\{x_1\in\Omega_1:\, \,
 (x_1, x_2)\in V^c \r\}\r)\, d\,\p_2(x_2)
  \geq \delta\, \p_2(W^c),
$$
which means $\p_2(W)\geq 1- \beta/\delta$. The choice $\delta = \sqrt{\beta}$ completes the proof.
\qed

The following lemma is obvious.

\begin{lemma}\label{dec}
Let $x_1, \ldots, x_N \in \Rn$, then
$$
 \sum_{i\ne j} \la x_i, x_j \ra  =  2^{2-N}
  \sum_{I\subset \{1, ..., N\} }\sum_{i\in I}
  \sum_{j\in I^c} \la x_i, x_j \ra.
  $$
\end{lemma}

\bigskip

\noindent
{\bf Proof of Theorem~\ref{A_k}. }
From Lemma \ref{dec} we have
$$
  \left| \left|\sum_{i=1}^N a_i X_i \r|^2 - \sum_{i=1}^N a_i^2  \left|
  X_i\r|^2 \r|   =  2^{2-N} \left|\sum _{I\subset\{1,2,..., N\}}
  \la \sum _{i\in I} a_i X_i, \sum _{j\in I^c}  a_j X_j \ra\right| .
$$
We deduce that
$$
 B_k^2\leq  2^{2-N} \sup_{a\in U_k}
 \sum _{I\subset\{1,2,..., N\}} Q(a, I, I^c)
  \leq
 2^{2-N}  \sum _{I\subset\{1,2,..., N\}} \sup_{a\in U_k}  Q(a, I, I^c)
 $$
 $$
 \leq  2^{2-N}  \sum _{I\subset\{1,2,..., N\}} Q_k(I).
$$
Let  $I\subset\{1,\dots,N\}$ be fixed. Proposition~\ref{estforq}  implies 
\begin{equation} \label{fixedi}
   \p\left(  Q_k(I)\leq M_0   \r) \geq 1-\beta ,
\end{equation}
where
$$
 M_0 := C_{\phi} \, t\, \max_{i\leq N} |X_i| + (M_1/4) \,  A_{k} .
$$

Consider two probability spaces $\{I : \, I\subset \{1, ..., N\}\}$ with the normalized
counting measure $\mu$ and our initial probability space $(\Omega, \p)$, on which $X_i$'s
are defined. By (\ref{fixedi}) we observe that the
$\mu\otimes \p$ probability of the event $V:=\{Q_k(I)\leq M_0  \}$ is at least $1-\beta$.
Then Lemma~\ref{fubini} implies that there exists $W\subset \Omega$ such that
$\p(W) \geq 1-\sqrt{\beta}$ and such that for every $\omega \in W$ one has
$\mu (\{Q_k(I) \leq M_0 \}) \geq 1-\sqrt{\beta}$.
 Since $Q_k(I) \leq A_k^2$, we obtain that for every $\omega \in W$,
$$
   B_k^2\leq  4 M_0 +4 \sqrt\beta A_k^2.
$$
Since $A_k^2 \leq \max_{i\leq N}  |X_i|^2+ B_k^2,$  we have
\begin{equation}\label{aandbk}
 A_k^2 \leq \frac{4M _0 + \max_{i\leq N}  |X_i|^2}{1-4
   \sqrt\beta}\quad \mbox{ and }
 \quad B_k^2 \leq \frac{4(M _0+ \sqrt{\beta} \max_{i\leq N}  |X_i|^2)}{1-4 \sqrt\beta}.
\end{equation}
Therefore
$$
  A_k^2 \leq (1-4\sqrt{\beta})^{-1} \left(\max_{i\leq N} |X_i|^2 + 4 C_{\phi}
  t \max_{i\leq N} |X_i| + M_1 A_k  \r)  .
$$
Using $\sqrt{u^2+v^2}\leq u+v$, and denoting
$\gamma =(1-4\sqrt{\beta})^{-1}$ (recall $M= \max_{i\leq N}|X_i|$) we obtain
$$
   A_k \leq \sqrt{\gamma} \, M   + 2 \sqrt{C_{\phi}\, \gamma\, t \, M} + \gamma \, M_1,
$$
which proves the estimate for $A_k$. Plugging this into (\ref{aandbk}),
we also observe
\begin{align*}
  B_k^2 &\leq \gamma \left(4 \sqrt{\beta}\, M^2 + 4 C_{\phi} t M + \gamma M_1^2 +
  \sqrt{\gamma} \, M\, M_1 +  2 \sqrt{C_{\phi}\, \gamma\, t  M }\, M_1    \r) \\
  &\leq \gamma \left(4 \sqrt{\beta}\, M^2 + 8 C_{\phi} t M + 2 \gamma M_1^2 +
  \sqrt{\gamma} \, M\, M_1  \r) .
\end{align*}
This completes the proof.
\qed

\section{Optimality}
\label{sec:optimal}

In this section we discuss optimality of estimates in 
Theorems~\ref{A_k} and \ref{thm:RIP-M}.
In Propositions~\ref{6.5}, \ref{6.6} and \ref{6.7}
we will prove results justifying remarks on optimality
following these theorems.

\medskip

To obtain the lower estimates on $A_m$ we use the following observation.

\begin{lemma}\label{optA_k}
Let $A=(X_{i j})_{{i\leq n},{j\leq N}}$ be an $n\times N$ matrix with i.i.d.
entries.  Then
\begin{equation}\label{equ:condition on t}
\p (A_m\geq t)\geq \frac{1}{2}\quad \text{ whenever}\quad
\p \left(|X_{11}|\geq  \frac{t}{\sqrt m}\right)\geq \frac{m+1}{N}.
 \end{equation}
\end{lemma}

\noindent{\bf Proof.}
 For every $i\leq N$, let $X_j\in\R^n$ be the $j$-th columns of $A$.
For $ m\leq N$ we have
$$
  A_m=\sup_{a\in U_m} \left|\sum_{j=1}^N a_j X_j\r|\geq \sup_{a\in U_m}
  \left|\sum_{j=1}^N a_j X_{1j}\r|\geq \sup_{a\in U_m \atop a_j\in \{\pm 1/\sqrt m, 0\}}
  \left|  \sum_{j=1}^N a_j X_{1j}\r|
$$
$$
 = \frac{1}{\sqrt m} \sum_{j=1}^m X_{1 j}^*\geq \sqrt m \, X_{1 m}^*.
$$
Therefore, using independence, we have
$$
   \p \left(A_m \geq t \right)\geq  \p \left(X_{1 m}^*\geq \frac{t}{\sqrt m} \right)=\p(Y\geq m),
$$
 where $Y$ is a real random variable with a binomial distribution
 of size $N$ and parameter $v=\p (|X_{11}|\geq  \frac{t}{\sqrt m})$.
 It is well known that the median of $Y$, med($Y$) satisfies
 $$\lfloor Nv\rfloor \leq \text{med }(Y)\leq \lceil Nv\rceil.$$
 Thus
 $\p (A_m\geq t)\geq \frac{1}{2}$ whenever
 $m\leq \lfloor Nv\rfloor.$ This implies the result.
\qed

To evaluate RIP, we will use the following simple observation.

\begin{lemma}\label{rip-sharp}
Let $n\leq N$ and $m\leq N$. Let $A$ be an $n\times N$ random matrix satisfying
$$
  \p (A_m\geq t \sqrt{m})\geq \frac{1}{2}.
$$
Assume also that $A$ satisfies RIP$_m(\delta)$ for some $\delta <1$
with probability greater than $1/2$. Then
$$
   m t^2 \leq 2 n .
$$
\end{lemma}

\noindent{\bf Proof.}
As $A$ satisfies RIP$_m(\delta)$ for some $\delta <1$
with probability greater than $1/2$, then clearly
$$
  A_m^2=\sup_{a\in U_m} |\sum a_iX_i|^2\leq 2n
$$
with probability greater than $1/2$. Therefore, with positive probability one has
$$
 t \sqrt{m}  \leq A_m \leq \sqrt{2n} ,
$$
which implies the result.
\qed

In order to show that a matrix with i.i.d. random variables satisfies
condition  ${\bf H}(\phi)$ with $\phi(t)=t^p$ we need the
Rosenthal's inequality (\cite{Ros}, see also \cite{FHJSZ}). As usual, by $\|\cdot\|_q$
for a random variable $\xi$ we mean its $L_q$-norm and for an $a\in \R^n$ its
$\ell _q$-norm, that is
$$
  \|\xi\|_q = \left(\E |\xi|^q\r)^{1/q} \quad \mbox{ and } \quad
 \|a\|_q = \left(\sum _{i=1}^{n} |a_i|^q\r)^{1/q} .
$$
Note that originally the Rosenthal inequality was proved for symmetric random
variables, but using standard symmetrization argument (i.e., passing from random
variables $\xi_i$'s to $(\xi_i-\xi_i')$'s, where $(\xi_i')$'s have the same distribution
and are independent), one can pass to centered random variables.

\begin{lemma}\label{rosin-l} Let $q>2$ and $a\in \R^n$. Let $\xi_1$, ..., $\xi_n$ be i.i.d. centered
random variables with finite $q$-th moment. Then there exists a positive absolute constant $C$
such that
\begin{equation}\label{rosin}
  \frac{1}{2} \, M_q \leq \left\|\sum _{i=1}^{n} a_i \xi_i \r\|_q \leq C \frac{q}{\ln q}\, M_q
\end{equation}
where $M_q:=\max\left\{\|a\|_2 \|\xi_1\|_2 , \|a\|_q \|\xi_1\|_q\r\}$.
\end{lemma}

The following is an almost immediate corollary of Rosenthal's inequality.
It should be compared with Proposition~1.3 of \cite{SV}.

\begin{cor}\label{ros-conc} Let $p>4$.
Let $\xi$ be a random variable of variance one and with a finite $p$-th moment.
Let $\xi_{ij}$, $i\leq n$, $j\leq N$ be i.i.d. random variables distributed as
$\xi$. Then for every $t>0$,
$$
 \p \left(\max_{j\leq N} \left|\frac{1}{n} \sum _{i=1}^n \xi _{ij}^2 - 1\r|> t \r)
 \leq \left(\frac{C p}{t\, \ln p} \r)^{p/2}\,  \E |\xi |^p \, \, \frac{N}{n^{p/4}},
$$
where $C$ is a positive absolute constant.
\end{cor}

\noindent{\bf Proof.}
Let $\xi_1$, ..., $\xi_n$ be i.i.d. random variables distributed as
$\xi$. We apply Rosenthal's inequality to random variables $(\xi _i^2 -1)$
with $q=p/2$ and $a=(1,1,...,1)$. Then
$$
 \left\|\sum _{i=1}^{n}  (\xi_i^2 -1) \r\|_{p/2} \leq C_p \sqrt{n} \, \|\xi^2 -1\|_{p/2}
 \leq   C_p \sqrt{n} \, \left(\|\xi^2\|_{p/2} +1\r) \leq   2 C_p \sqrt{n} \, \|\xi\|_{p}^2,
$$
where $C_p=Cp/\ln p$ for an absolute positive constant $C$.
Using Chebyshev's inequality we observe
$$
  \p \left( \left|\frac{1}{n} \sum _{i=1}^n \xi _{i}^2 - 1\r|> t \r)
 \leq \frac{\E \sum_{i=1}^{n}|\xi_i^2 -1|^{p/2}}{(t n )^{p/2}} \leq
  \frac{(2 C_p)^{p/2}\, \|\xi\|_{p}^p}{t^{p/2}\, n^{p/4}} .
$$
The result follows by the union bound.
\qed

As is mentioned in remarks on optimality following Theorem~\ref{A_k} 
the next proposition gives a lower bound for $A_m$ to be compared with 
Case~1 of Theorem~\ref{A_k}.

\begin{prop}\label{optimal A_m}\label{6.5}
Let $p>2$, $1\leq m\leq N$. There exists a sequence of  independent random  vectors $X_1,\cdots,X_N$
in $\R^n$ satisfying
\begin{equation}\label{4thmomentc}
  \quad \forall 1\leq i\leq N\  \forall
  a \in S^{n-1}\quad \E |\la X_i, a\ra |^p \leq 1
\end{equation}
and such that
$$
  \p \left(A_m\geq \frac{C p}{\ln p}\,  \sqrt m  \left(\frac{N}{m}\right)^{1/p}
  \left(\ln\left(\frac{2 N}{m}\right)\right)^{-1/p}\right)\geq \frac{1}{2},
$$
  where $C$ is an absolute positive constant.
\end{prop}

\noindent{\bf Proof.}
Let $\lambda\geq 1$ to be set later and let us put
$$
 f_p(x)=
 \begin{cases}
 \frac{p}{2(1-\lambda^{-p})|x|^{p+1}} &\mbox{if } 1\leq |x|\leq \lambda \\
 0 &\mbox{otherwise}.
 \end{cases}
$$
We have
$\int f_p (x) \, dx=1$ and
$$
  a_p^p := \int |x|^pf_p(x)\, dx= p \frac{\ln \lambda}{1-\lambda^{-p}}.
$$
Consider the random variable $\xi (\omega) = \omega$ with respect to the density $f_p$ 
 and let $(X_{ij})$ be i.i.d. copies of  $\xi/a_p$.
Clearly, $\E |X_{11}|^p=1$. Since, for $s\in[1,\lambda]$
$$
 \p \left(|\xi|>s\right)= \frac{1}{1-\lambda^{-p}}\left(\frac{1}{s^p}-\frac{1}{\lambda^p}\right),
$$
a short computation using (\ref{equ:condition on t}) shows that $\p (A_m\geq t)\geq \frac{1}{2}$
provided that
$$
  t\leq  \left(\frac{1-\lam^{-p}}{p \, \ln \lam} \r)^{1/p}  \sqrt m
  \left(\frac{N}{(m+1)(1-\lam^{-p}) + N \lam^{-p}}\right)^{1/p}
$$
$$
 =
 \sqrt m\, \left(\frac{1}{p \, \ln \lam} \r)^{1/p}
  \left(\frac{N}{m+1 + N/(\lam^{p}-1)}\right)^{1/p}.
$$
Choosing $\lam$ from $\lam ^p -1 = N/(m+1)$, we obtain $\p (A_m\geq t)\geq \frac{1}{2}$
provided that
$$
 t \leq \sqrt{m} \left(\frac{N}{2(m+1) \ln (2N/(m+1))}\r)^{1/p}.
$$
Finally, to satisfy condition (\ref{4thmomentc}), we pass from matrix $A$ to
$A'=A/c_p = ( X_{ij}/c_p ) _{ij}$, where $c_p\leq C p/\ln p$ is a constant
in Rosenthal's inequality (\ref{rosin}).
By Rosenthal's inequality, the sequence of columns of $A'$ satisfies the condition
 (\ref{4thmomentc}).
\qed

The next proposition gives an upper bound on the size of sparsity $m$ in order to satisfy RIP under condition of Case 1 of Theorem~\ref{thm:RIP-M}  (see  Remark~3 following
this theorem).

\begin{prop}\label{opt-rip-4m}\label{6.6}
Let $q>p>2$, $n\leq N$ and $m\leq N$. There exist an absolute positive constant  $C$,
an $n\times N$ matrix $A$, whose columns $X_1, ..., X_N$ are  independent random vectors
satisfying
\begin{equation}\label{pmomen}
  \quad \forall 1\leq i\leq N\  \forall
  a \in S^{n-1}\quad \E |\la X_i, a\ra |^p \leq \left(\frac{C p}{\ln p}\r)^p
  \, \frac{q}{q-p}\, \left( \frac{q-2}{q}\r)^{p/2} ,
\end{equation}
and  for every $t\in (0, 1)$,
\begin{equation}\label{rip-conc-cond}
 \p\left(\max_{i\leq N} \left|\frac{|X_i|^2}{n}-1 \r|\geq t\right)\leq t^{p/2}
\end{equation}
provided that
$$
  N \leq \left(\frac{q\ln p}{C (q-2) p}\r)^{p/2} \, \frac{q-p}{q}\, t^{p} \, n^{p/4}.
$$
Assume that  $A$ satisfies RIP$_m(\delta)$ for some $\delta <1$
with probability greater than $1/2$.  Then
$$
   m \, \left(\frac{N}{m+1}\r)^{2/q}\leq  \frac{2(q-2)}{q}\,  n .
$$
\end{prop}

\noindent{\bf Proof.}
Consider the density
$$
 f(x)=
 \begin{cases}
 \frac{q}{2 |x|^{q+1}} &\mbox{if }  |x|\geq 1 \\
 0 &\mbox{otherwise} .
 \end{cases}
$$
We have
$\int f (x) \, dx=1$,
$$
  \int |x|^p f(x)\, dx=  \frac{q}{q-p} \quad \mbox{ and } \quad
  a_2^2:= \int |x|^2 f(x)\, dx=  \frac{q}{q-2}.
$$
Consider the random variable $\xi (\omega) = \omega$ with respect to the
density $f$ and let $(X_{ij})_{ij}$ be i.i.d. copies of  $\xi/a_2$.
Clearly,
$$
  \E |X_{11}|^2=1 \quad  \mbox{ and } \quad
  \E |X_{11}|^p=\frac{q}{q-p}\, \left(\frac{q-2}{q}\r)^{p/2}.
$$
Then Rosenthal's inequality (\ref{rosin}) implies the condition (\ref{pmomen})
and Corollary~\ref{ros-conc} implies (\ref{rip-conc-cond}).

Now we estimate $A_m$ for the matrix $A$, whose columns are $(X_{ij})_i$, $j\leq N$.
Since, for $s\geq 1$, $\p \left(|\xi|>s\right)= s^{-q}$, by
(\ref{equ:condition on t}), we obtain  that $\p (A_m\geq t)\geq \frac{1}{2}$
provided that
$$
  t\leq   \sqrt m \, \sqrt{\frac{q-2}{q}} \, \left(\frac{N}{m+1} \r)^{1/q}.
$$
This means
$$
   \p \left(A_m\geq \sqrt m \, \sqrt{\frac{q-2}{q}} \, \left(\frac{N}{m+1} \r)^{1/q}\r)
   \geq \frac{1}{2},
$$
and we complete the proof applying Lemma~\ref{rip-sharp}.
\qed

The next proposition shows the optimality (up to absolute constants)
of the sparsity parameter in Case~2 of Theorem~\ref{thm:RIP-M}
(see Remark~4 following this theorem) as well as optimality of bounds for $A_m$ in Case~2 of Theorem~\ref{A_k} 
(see remarks on optimality following this theorem).

\begin{prop}\label{optimal A} \label{6.7}
There exist absolute positive constants $c$, $C$
such that the following holds.
Let $\alpha \in [1,2]$, $1\leq m\leq N/2$ and $n$ satisfies $N\leq \exp(cn^{\alpha/2})$.
There exists an $n\times N$ matrix $A$, whose columns $X_1, ..., X_N$ are  independent random vectors satisfying
\begin{equation}\label{condonepsi}
  \quad \forall 1\leq i\leq N\, \, \,  \forall
  a \in S^{n-1}\quad \E \exp\left(|\la X_i, a\ra |^{\alpha}\r) \leq C
\end{equation}
and
\begin{equation}\label{condtwoconc}
 \p\left(\max_{i\leq N} \left|\frac{|X_i|^2}{n}-1 \r|\geq \frac{\sqrt{2}-1}{2}\right)
  \leq 2 \exp(-c n^{\alpha/2}),
\end{equation}
and such that
\begin{equation}\label{sharpest}
  \p \left(A_m\geq \sqrt{\frac{m}{2}}\, \left(
  \ln \frac{N}{m+1} \right)^{1/\alpha} \right)\geq \frac{1}{2}.
\end{equation}
Additionally, if $n\leq N$ and if $A$ satisfies RIP$_m(\delta)$ for some $\delta <1$
with probability greater than $1/2$, then
$$
   m \left( \ln \frac{N}{m+1}\r)^{2/\alpha}\leq 4 n .
$$
\end{prop}

\noindent{\bf Proof.}
We consider a symmetric random variable $\xi$ with the distribution defined by
$\p\left(|\xi |> t\r)= \exp(-t^{\alpha})$. It is easy to check that
$$
  \E \exp(|\xi|^{\alpha} / 2) =2
$$
and
$$
   a:= \E \xi^2 =  \Gamma\left(\frac{2}{\alpha}+1\r)\in [1, 2].
$$
Let $X_{ij}$, $i\leq n$, $j\leq N$ be i.i.d. copies of $\xi/ \sqrt{a}$, $A=(X_{ij})_{ij}$ and
$X_j$'s be its columns.
Applying Lemma~3.4 from \cite{ALPT2} (see also Theorem~1.2.8 in \cite{lama})
we observe that $X_i$'s satisfy conditions
(\ref{condonepsi}) and (\ref{condtwoconc}).
By (\ref{equ:condition on t}) we observe that $\p (A_m\geq t)\geq \frac{1}{2}$
provided that
$$
   \p\left(|\xi|\geq \frac{\sqrt{a}\, t}{\sqrt{m}}\r) =
   \exp\left(-\left(\sqrt{a} t/\sqrt{m}\r)^{\alpha}\r) \geq \frac{m+1}{N}.
$$
Thus it is enough to take
$$
   t\leq \sqrt{\frac{{m}}{a}}\, \left( \ln \frac{N}{m+1} \r)^{1/\alpha} .
$$
This proves the estimate (\ref{sharpest}).

Finally, the ``additionally" part follows by Lemma~\ref{rip-sharp}.
\qed

\address
\end{document}